\documentclass[11pt]{article}
\usepackage{amsmath,amssymb,amsfonts, amsbsy, amsthm, epsfig, subfig, graphicx,rotating}
\usepackage[usenames]{color}

\newtheorem{theorem}{Theorem}
\newtheorem{corollary}{Corollary}
\newtheorem{proposition}{Proposition}
\newtheorem{lemma}{Lemma}
{
\theoremstyle{definition}
\newtheorem{definition}{Definition}
\newtheorem{example}{Example}
\newtheorem{remark}{Remark}
}
\newcommand{\beq}{\begin{equation}}
\newcommand{\eeq}{\end{equation}}
\newcommand{\beas}{\begin{eqnarray*}}
\newcommand{\eeas}{\end{eqnarray*}}
\newcommand{\bea}{\begin{eqnarray}}
\newcommand{\eea}{\end{eqnarray}}
\newcommand{\bei}{\begin{itemize}}
\newcommand{\eei}{\end{itemize}}
\newcommand{\ben}{\begin{enumerate}}
\newcommand{\een}{\end{enumerate}}
\newcommand{\bet}{\begin{theorem}}
\newcommand{\eet}{\end{theorem}}
\newcommand{\bel}{\begin{lemma}}
\newcommand{\eel}{\end{lemma}}
\newcommand{\bep}{\begin{proposition}}
\newcommand{\eep}{\end{proposition}}
\newcommand{\bed}{\begin{definition}}
\newcommand{\eed}{\end{definition}}
\newcommand{\bec}{\begin{corollary}}
\newcommand{\eec}{\end{corollary}}
\newcommand{\bex}{\begin{example}}
\newcommand{\eex}{\end{example}}

\newcommand{\argmin}{\mathop{\rm arg\min}}

\begin{document}

\title{Estimating Sparse Precision Matrix: Optimal Rates of Convergence and Adaptive Estimation}
\author{T. Tony Cai$^{1},$ Weidong Liu$^{2}$ and Harrison H. Zhou$^{3}$}
\date{}
\maketitle

\begin{abstract}

Precision matrix is of significant importance in a wide range of applications in multivariate analysis.
This paper considers adaptive minimax estimation of sparse precision matrices in the
high dimensional setting. Optimal rates of convergence are established for a range of matrix norm losses. A
fully data driven estimator based on adaptive constrained $\ell_1$
minimization is proposed and its rate of convergence is obtained over
a collection of parameter spaces. The estimator, called ACLIME, is easy to
implement and performs well numerically.

A major step in establishing the minimax rate of convergence is the derivation of a rate-sharp lower bound.
A ``two-directional" lower bound technique is applied
to obtain the minimax lower bound. The upper and lower bounds together yield the optimal rates of
convergence  for sparse precision matrix estimation and show that the ACLIME estimator is adaptively minimax rate optimal for a collection of parameter spaces and a range of matrix norm losses simultaneously.
\end{abstract}

\footnotetext[1]{%
Department of Statistics, The Wharton School, University of Pennsylvania,
Philadelphia, PA 19104. \newline
\indent\ \ The research of Tony Cai was supported in part by NSF Grant
DMS-0604954 and NSF FRG Grant \newline
\indent\ \ DMS-0854973.} \footnotetext[2]{%
Department of Mathematics and Institute of Natural Sciences, Shanghai Jiao
Tong University, Shanghai, \newline
\indent\ \ China.} \footnotetext[3]{%
Department of Statistics, Yale University, New Haven, CT 06511. The research
of Harrison Zhou was \newline
\indent\ \ supported in part by NSF Career Award DMS-0645676 and NSF
FRGGrant DMS-0854975.}



\newpage

\section{Introduction}
\label{intro.sec}

Precision matrix plays a fundamental role in many high-dimensional inference
problems. For example, knowledge of the precision matrix is crucial for
classification and discriminant analyses. Furthermore, precision matrix
is critically useful for a broad range of applications such as portfolio
optimization, speech recognition, and genomics. See, for example, Lauritzen
(1996), Yuan and Lin (2007), Saon and Chien (2011).
Precision matrix is also closely connected to the graphical models
which are a powerful tool to model the relationships among a large number of
random variables in a complex system and are used in a wide array of
scientific applications. It is well known that recovering the structure of
an undirected Gaussian graph is equivalent to the recovery of the support of the precision matrix. See for example, Lauritzen (1996), Meinshausen and B\"{u}hlmann (2006) and Cai, Liu and Luo (2011). Liu, Lafferty and Wasserman (2009) extended the result to a more general class of distributions called nonparanormal distributions.

The problem of estimating a large precision matrix and recovering its
support has drawn considerable recent attention and a number of methods have
been introduced. Meinshausen and B\"{u}hlmann (2006) proposed a neighborhood
selection method for recovering the support of a precision matrix. Penalized
likelihood methods have also been introduced for estimating sparse precision
matrices. Yuan and Lin (2007) proposed an $\ell _{1}$ penalized normal
likelihood estimator and studied its theoretical properties. See also
 Friedman, Hastie and Tibshirani (2008), d'Aspremont, Banerjee and El Ghaoui (2008), Rothman et al. (2008),
Lam and Fan (2009), and Ravikumar et al. (2011). Yuan (2010) applied the
Dantzig Selector method to estimate the precision matrix and gave the
convergence rates for the estimator under the matrix $\ell _{1}$ norm and
spectral norm. Cai, Liu and Luo (2011) introduced an estimator called CLIME
using a constrained $\ell _{1}$ minimization approach and obtained the rates
of convergence for estimating the precision matrix under the spectral norm
and Frobenius norm.

Although many methods have been proposed and various rates of convergence
have been obtained, it is unclear which estimator is optimal for estimating a
sparse precision matrix in terms of convergence rate. This is due to the
fact that the minimax rates of convergence, which can serve as a fundamental
benchmark for the evaluation of the performance of different procedures, is
still unknown. The goals of the present paper are to establish the optimal
minimax rates of convergence for estimating a sparse precision matrix under
a class of matrix norm losses and to introduce a fully data driven adaptive
estimator that is simultaneously rate optimal over a  collection of parameter spaces for each loss in this class.

Let $X_{1},\ldots ,X_{n}$ be a random sample from a $p$-variate distribution
with a covariance matrix $\Sigma =\left( \sigma _{ij}\right) _{1\leq i,j\leq
p}$. The goal is to estimate the inverse of $\Sigma$, the precision
matrix $\Omega =\left( \omega _{ij}\right) _{1\leq i,j\leq p}$. It is well
known that in the high-dimensional setting structural assumptions are needed
in order to consistently estimate the precision matrix. The class of sparse
precision matrices, where most of the entries in each row/column are zero or
negligible, is of particular importance as it is related to sparse graphs in the Gaussian case.
For a matrix $A$ and
a number $1\leq r\leq \infty $, the matrix $\ell _{w}$ norm is defined as $%
\Vert A\Vert _{w}=\sup_{|x|_{w}\leq 1}|Ax|_{w}$. In particular, the commonly
used spectral norm is the matrix $\ell _{2}$ norm. For a symmetric matrix $A$%
, it is known that the spectral norm $\left\Vert A\right\Vert _{2}$ is equal
to the largest magnitude of eigenvalues of $A$.
The sparsity of a precision matrix
can be modeled by the $\ell _{q}$ balls with $0\leq q<1$.
More specifically,
we define the parameter space $\mathcal{G}_{q}(c_{n,p},M_{n,p})$ by
\begin{equation}
\mathcal{G}_{q}(c_{n,p},M_{n,p})=\left\{
\begin{array}{c}
\Omega =\left( \omega _{ij}\right) _{1\leq i,j\leq
p}:\max_{j}\sum_{i=1}^{p}|\omega _{ij}|^{q}\leq c_{n,p}, \\
\Vert \Omega \Vert _{1}\leq M_{n,p},\text{ }\lambda _{\max }(\Omega
)/\lambda _{\min }(\Omega )\leq M_{1},\Omega \succ 0%
\end{array}%
\right\} ,  \label{parameterspace}
\end{equation}%
where $0\leq q<1$, $M_{n,p}$ and $c_{n,p}$ are positive and bounded away
from $0$, $M_1>0$ is a given constant,  $\lambda _{\max }(\Omega)$ and $\lambda _{\min }(\Omega)$ are the largest and smallest eigenvalues of $\Omega$ respectively, and $c_{1}n^{\beta }\leq p\leq \exp \left( \gamma n\right) $ for some
constants $\beta >1$, $c_{1}>0$ and $\gamma >0$. The notation $A\succ 0$ means that $A$
is symmetric and positive definite. In the special case of $q=0$, a matrix
in $\mathcal{G}_{0}(c_{n,p},M_{n,p})$ has at most $c_{n,p}$ nonzero elements
on each row/column.

Our analysis establishes the minimax rates of convergence for estimating the
precision matrices over the parameter space $\mathcal{G}_{q}(c_{n,p}, M_{n,p})$ under the matrix $\ell _{w}$ norm losses for $1\le w \le \infty$. We shall first introduce a new method using an adaptive
constrained $\ell_1$ minimization approach for estimating the sparse
precision matrices. The estimator, called ACLIME, is fully data-driven and
easy to implement. The properties of the ACLIME are then studied in detail
under  the matrix $\ell _{w}$ norm losses. In particular, we establish the rates
of convergence for the ACLIME estimator  which provide
upper bounds for the minimax risks.

A major step in establishing the minimax rates of convergence is the
derivation of rate sharp lower bounds. As in the case of estimating sparse
covariance matrices, conventional lower bound techniques, which are designed
and well suited for problems with parameters that are scalar or
vector-valued, fail to yield good results for estimating sparse precision
matrices under the spectral norm. In the present paper we apply the
``two-directional" lower bound technique first developed in Cai and Zhou
(2012) for estimating sparse covariance matrices. This lower bound method
can be viewed as a simultaneous application of Assouad's Lemma along the row
direction and Le Cam's method along the column direction. The lower bounds
match the rates in the upper bounds for the ACLIME estimator and thus yield
the minimax rates.

By combining the minimax lower and upper bounds developed in later sections,
the main results on the optimal rates of convergence for estimating a sparse
precision matrix under various norms can be summarized in the following
theorem. We focus here on the exact sparse case of $q=0$; the optimal rates for the general case of $0\le q <1$ are given in the end of Section \ref{lowerbd.sec}. Here for two sequences of positive numbers $a_{n}$ and $b_{n}$, $%
a_{n}\asymp b_{n}$ means that there exist positive constants $c$ and $C$
independent of $n$ such that $c\leq a_{n}/b_{n}\leq C$.

\begin{theorem}
\label{MinimaxOpe.q0}
Let $X_{i}{\overset{\text{i.i.d.}}{\sim }}N_{p}(\mu ,\Sigma )$, $i=1,2,\ldots ,n,$ and let
$1 \leq k =o(n^{\frac{1}{2}}\left( \log p\right) ^{-\frac{3}{2}}).$
The minimax risk of estimating the precision matrix $\Omega =\Sigma^{-1}$  over the class
$\mathcal{G}_{0}(k, M_{n,p})$ based on the random sample $\{X_1, ..., X_n\}$ satisfies
\begin{equation}
\inf_{\hat{\Omega}}\sup_{\mathcal{G}_{0}(k,M_{n,p})}\mathbb{E}%
\left\Vert \hat{\Omega}-\Omega \right\Vert _{w}^{2}\asymp
M_{n,p}^{2} k^{2} \frac{\log p}{n}
\label{rateOper}
\end{equation}%
for all $1\leq w\leq \infty$.
\end{theorem}

In view of Theorem \ref{MinimaxOpe.q0}, the ACLIME estimator, which is fully
data driven, attains the optimal rates of convergence simultaneously for all $k$-sparse precision matrices in
 the parameter spaces $\mathcal{G}_{0}(k, M_{n,p})$ with $k \ll n^{\frac{1}{2}}\left( \log p\right) ^{-\frac{3}{2}}$ under the matrix $\ell_w$ norm for all $1\le w \le \infty$. As will be seen in Section \ref{lowerbd.sec}, the adaptivity holds for the  general  $\ell_q$ balls $\mathcal{G}_{q}(c_{n,p}, M_{n,p})$ with $0\le q< 1$.
The ACLIME procedure is thus  rate optimally adaptive to both the sparsity
patterns and the loss functions.

In addition to its theoretical optimality, the ACLIME estimator is
computationally easy to implement for high dimensional data. It can be
computed column by column  via linear programming and the algorithm
is easily scalable. A simulation study is carried out to investigate the numerical performance of the ACLIME estimator.  The results show that the procedure performs favorably in comparison to CLIME.





Our work on optimal estimation of precision matrix given in the present
paper is closely connected to a growing literature on estimation of large
covariance matrices. Many regularization methods have been proposed and
studied. For example, Bickel and Levina (2008a, b) proposed banding and
thresholding estimators for estimating bandable and sparse covariance
matrices respectively and obtained rate of convergence for the two
estimators. See also El Karoui (2008) and Lam and Fan (2009). Cai, Zhang and
Zhou (2010) established the optimal rates of convergence for estimating
bandable covariance matrices. Cai and Yuan (2012) introduced an adaptive
block thresholding estimator which is simultaneously rate optimal rate over
large collections of bandable covariance matrices. Cai and Zhou (2012)
obtained the minimax rate of convergence for estimating sparse covariance
matrices under a range of losses including the spectral norm loss. In
particular, a new general lower bound technique was developed. Cai and Liu
(2011) introduced an adaptive thresholding procedure for estimating sparse
covariance matrices that automatically adjusts to the variability of
individual entries.

The rest of the paper is organized as follows.  The ACLIME estimator is introduced in detail in Section \ref{method.sec}  and its theoretical properties are studied in Section \ref{optimality.sec}. In particular, a minimax upper bound for estimating sparse precision matrices is obtained.
Section \ref{lowerbd.sec}
establishes a minimax lower bound which matches the minimax upper bound derived in Section \ref{method.sec} in terms of the convergence rate. The
upper and lower bounds together yield the optimal minimax rate of
convergence. A simulation study is carried out
in Section \ref{simulation.sec} to
compare the performance of the ACLIME with that of the CLIME estimator. Section \ref{discussion.sec} gives the optimal rate of convergence for estimating sparse precision matrices under the Frobenius norm and discusses connections and differences
of our work with other related problems.  The proofs are given in Section \ref{proof.sec}.


\section{Methodology}
\label{method.sec}


In this section we introduce an adaptive constrained $\ell_1$ minimization
procedure, called ACLIME, for estimating a precision matrix $\Omega$. The properties of the estimator are then studied in Section \ref{optimality.sec} under the matrix $\ell_w$ norm losses for $1\le w\le \infty$  and a minimax upper bound is established.
The upper bound together with the lower bound given in Section \ref{lowerbd.sec} will show that the
ACLIME estimator is adaptively rate optimal.

We begin with basic notation and definitions. For a vector $a=(a_{1},\dotsc
,a_{p})^{T}\in \mathbb{R}^{p}$, define $|a|_{1}=\sum_{j=1}^{p}|a_{j}|$ and $%
|a|_{2}=\sqrt{\sum_{j=1}^{p}a_{j}^{2}}$. For a matrix $A=(a_{ij})\in \mathbb{%
R}^{p\times q}$, we define the elementwise $\ell _{r}$ norm by $%
|A|_{r}=(\sum_{i,j}|a_{ij}|^{r})^{1/r}$. The Frobenius norm of $A$ is the
elementwise $\ell _{2}$ norm. $I$ denotes a $p\times p$ identity matrix. For
any two index sets $T$ and $T^{^{\prime }}$ and matrix $A$, we use $%
A_{TT^{^{\prime }}}$ to denote the $|T|\times |T^{^{\prime }}|$ matrix with
rows and columns of $A$ indexed by $T$ and $T^{^{\prime }}$ respectively.




For an i.i.d. random sample $\{X_{1},\ldots ,X_{n}\}$ of $p$-variate
observations drawn from a population $X$, let the sample mean $\bar{X}=\frac{1}{n}\sum_{k=1}^{n}X_{k}$ and the sample covariance matrix
\begin{equation}
\Sigma ^{\ast }=(\sigma _{ij}^{\ast })_{1\leq i,j\leq p}=\frac{1}{n-1}%
\sum_{l=1}^{n}\left( X_{l}-\bar{X}\right) \left( X_{l}-\bar{X}\right) ^{T},
\label{MLE}
\end{equation}%
which is an unbiased estimate of the covariance matrix $\Sigma =\left(
\sigma _{ij}\right) _{1\leq i,j\leq p}$. 

It is well known that in the high dimensional setting, the inverse of the
sample covariance matrix either does not exist or is not a good estimator of
$\Omega $. As mentioned in the introduction, a number of methods for
estimating $\Omega $ have been introduced in the literature. In particular,
Cai, Liu and Luo (2011) proposed an estimator called CLIME by solving the
following optimization problem:
\begin{equation}
\min |\Omega|_{1}~~\mbox{subject to:}~~|\Sigma ^{\ast }\Omega
-I|_{\infty }\leq \tau _{n},~~\Omega \in \mathbb{R}^{p\times p},  \label{c1}
\end{equation}%
where $\tau _{n}=CM_{n,p}\sqrt{\log p/n}$ for some constant $C$. The convex
program (\ref{c1}) can be further decomposed into $p$ vector-minimization
problems. Let $e_{i}$ be a standard unit vector in $\mathbb{R}^{p}$ with $1$
in the $i$-th coordinate and $0$ in all other coordinates. For $1\leq i\leq
p $, let $\hat{\omega}_{i}$ be the solution of the following convex
optimization problem
\begin{equation}
\min |\omega |_{1}~~\mbox{subject to}~~|\Sigma _{n}\omega -e_{i}|_{\infty
}\leq \tau _{n},  \label{o1}
\end{equation}%
where $\omega $ is a vector in $\mathbb{R}^{p}$. The final CLIME estimator
of $\Omega $ is obtained by putting the columns $\hat{\omega}_{i}$ together
and applying an additional symmetrization step. This estimator is easy to
implement and possesses a number of desirable properties as shown in Cai,
Liu and Luo (2011).

The CLIME estimator has, however, two drawbacks. One is that the estimator
is not rate optimal, as will be shown later. Another drawback is that the
procedure is not adaptive in the sense that the tuning parameter $\lambda
_{n}$ is not fully specified and needs to be chosen through  an empirical
method such as cross-validation.


To overcome these drawbacks of CLIME, we now introduce an adaptive
constrained $\ell_1$-minimization for inverse matrix estimation (ACLIME).
The estimator is fully data-driven and adaptive to the variability of
individual entries.
A key technical result which provides the motivation for the new procedure
is the following fact.

\begin{lemma}
\label{var.bnd.lem} Let $X_{1},...,X_{n}\overset{iid}{\sim }N_{p}(\mu
,\Sigma )$ with $\log p=O(n^{1/3})$. Set $S^{\ast }=(s_{ij}^{\ast })_{1\leq
i,j\leq p}=\Sigma ^{\ast }\Omega -I_{p\times p}$, where $\Sigma ^{\ast }$ is
the sample covariance matrix defined in (\ref{MLE}). Then%
\begin{equation*}
Var\left( s_{ij}^{\ast }\right) =\left\{
\begin{array}{cc}
n^{-1}(1+\sigma _{ii}\omega _{ii}), & \text{for }i=j \\
n^{-1}\sigma _{ii}\omega _{jj}, & \text{for }i\neq j%
\end{array}%
\right.
\end{equation*}%
and for all $\delta \geq 2$,
\begin{equation}
\mathbb{P}\left\{ |(\Sigma ^{\ast }\Omega -I_{p\times p})_{ij}|\leq \delta
\sqrt{\frac{\sigma _{ii}\omega _{jj}\log p}{n}},\;\forall \;1\leq i,j\leq
p\right\} \geq 1-O((\log p)^{-{\frac{1}{2}}}p^{-{\frac{\delta ^{2}}{4}}+1}).
\label{le1}
\end{equation}
\end{lemma}

A major step in the construction of the adaptive data-driven procedure is to
make the constraint in \eqref{c1} and \eqref{o1} adaptive to the variability
of individual entries based on Lemma \ref{var.bnd.lem}, instead of using a
single upper bound $\lambda _{n}$ for all the entries. In order to apply
Lemma \ref{var.bnd.lem}, we need to estimate the diagonal elements of $%
\Sigma $ and $\Omega $, $\sigma _{ii}$ and $w_{jj}$, $i,j=1,...,p$. Note
that $\sigma _{ii}$ can be easily estimated by the sample variances $\sigma
_{ii}^{\ast }$, but $\omega _{jj}$ are harder to estimate. Hereafter, $%
(A)_{ij}$ denotes the $(i,j)$-th entry of the matrix $A$, $(a)_{j}$ denotes
the $j$-th element of the vector $a$. Denote $b_{j}=(b_{1j},\ldots
,b_{pj})^{^{\prime }}$.

The ACLIME procedure has two steps: The first step
is to estimate $\omega _{jj}$ and the second step is to apply a modified
version of the CLIME procedure to take into account of the variability of
individual entries.

\begin{enumerate}
\item[\textbf{Step 1:}] \textbf{Estimating }$\omega _{jj}$. Note that $%
\sigma _{ii}\omega _{jj}\leq (\sigma _{ii}\vee \sigma _{jj})\omega _{jj}$
and $(\sigma _{ii}\vee \sigma _{jj})\omega _{jj}\geq 1$. So the inequality
on the left hand side of (\ref{le1}) can be relaxed to
\begin{equation}
|(\Sigma ^{\ast }\Omega -I_{p\times p})_{ij}|\leq 2(\sigma _{ii}\vee \sigma
_{jj})\omega _{jj}\sqrt{\frac{\log p}{n}},\quad 1\leq i,j\leq p.  \label{le2}
\end{equation}%
Let $\hat{\Omega}_{1}:=(\hat{\omega}_{ij}^{1})=(\hat{\omega}_{\cdot
1}^{1},\ldots ,\hat{\omega}_{\cdot p}^{1})$ be a solution to the following
optimization problem:
\begin{equation}
\hat{\omega}_{\cdot j}^{1}=\mathop{\rm arg\min}_{b_{j}\in R^{p}}\left\{
|b_{j}|_{1}:\;|\hat{\Sigma}b_{j}-e_{j}|_{\infty }\leq \lambda _{n}(\sigma
_{ii}^{\ast }\vee \sigma _{jj}^{\ast })\times b_{jj},\quad b_{jj}>0\right\} ,
\label{le3}
\end{equation}%
where $b_{j}=(b_{1j},\ldots ,b_{pj})^{^{\prime }}$, $1\leq j\leq p$, $\hat{%
\Sigma}=\Sigma ^{\ast }+n^{-1}I_{p\times p}$ and
\begin{equation}
\lambda _{n}=\delta \sqrt{\frac{\log p}{n}}.  \label{lambdan}
\end{equation}%
Here $\delta $ is a constant which can be taken as 2. The estimator $\hat{%
\Omega}_{1}$ yields estimates of the conditional variance $\omega _{jj}$, $%
1\leq j\leq p$. More specifically, we define the estimates of $\omega _{jj}$
by
\begin{equation*}
\breve{\omega}_{jj}=\hat{\omega}_{jj}^{1}I\left\{ \sigma _{jj}^{\ast }\leq
\sqrt{\frac{n}{\log p}}\right\} +\sqrt{\frac{\log p}{n}}I\left\{
\sigma _{jj}^{\ast }>\sqrt{\frac{n}{\log p}}\right\} .
\end{equation*}

\smallskip

\item[\textbf{Step 2:}] \textbf{Adaptive estimation. } Given the estimates $%
\breve{\omega}_{jj}$, the final estimator $\hat{\Omega}$ of $\Omega $ is
constructed as follows. First we obtain $\tilde{\Omega}^{1}=:(\tilde{\omega}%
_{ij}^{1})$ by solving $p$ optimization problems: for $1\leq j\leq p$
\begin{equation}
\tilde{\omega}_{\cdot j}^{1}=\mathop{\rm arg\min}_{b\in R^{p}}\left\{
|b|_{1}:\;|(\hat{\Sigma}b-e_{j})_{i}|\leq \lambda _{n}\sqrt{\sigma
_{ii}^{\ast }\breve{\omega}_{jj}},\quad 1\leq i\leq p\right\},  \label{op1}
\end{equation}%
where $\lambda _{n}$ is given in \eqref{lambdan}. We then obtain the
estimator $\hat{\Omega}$ by symmetrizing $\tilde{\Omega}^{1}$,
\begin{equation}
\hat{\Omega}=(\hat{\omega}_{ij}),~\text{where }\hat{\omega}_{ij}=\hat{%
\omega}_{ji}=\tilde{\omega}_{ij}^{1}I\{|\tilde{\omega}_{ij}^{1}|\leq |\tilde{%
\omega}_{ji}^{1}|\}+\tilde{\omega}_{ji}^{1}I\{|\tilde{\omega}_{ij}^{1}|>|%
\tilde{\omega}_{ji}^{1}|\}.  \label{procedure}
\end{equation}
\end{enumerate}

We shall call the estimator $\hat{\Omega}$ adaptive CLIME, or ACLIME. The
estimator adapts to the variability of individual entries by using an
entry-dependent threshold for each individual $\omega _{ij}$. Note that the
optimization problem (\ref{le3}) is convex and can be cast as a linear
program. The constant $\delta$ in \eqref{lambdan} can be taken as 2 and the
resulting estimator will be shown to be adaptively minimax rate optimal for
estimating sparse precision matrices.

\begin{remark}
Note that $\delta=2$ used in the constraint sets is tight, it can not be
further reduced in general. If one chooses the constant $\delta<2$, then with probability tending to
1, the true precision matrix will no longer belong to the feasible sets. To
see this, consider $\Sigma=\Omega=I_{p\times p}$ for simplicity. It follows
from Liu, Lin and Shao (2008) and Cai and Jiang (2011) that
\[
\sqrt{\frac{n}{\log p}}\max_{1\leq i<j\leq p}|\hat{\sigma}_{ij}|\rightarrow 2
\]
in probability. Thus $P(|\hat{\Sigma}\Omega-I_{p\times p}|_{\infty}>\lambda_{n})\rightarrow 1$, which means that if $\delta<2$, the
true $\Omega$ lies outside of the feasible set with high probability and
solving the corresponding minimization problem cannot lead to a good
estimator of $\Omega$.
\end{remark}

\begin{remark}
The CLIME estimator uses a universal tuning parameter $\lambda _{n}=CM_{n,p}%
\sqrt{\log p/n}$ which does not take into account the variations in the
variances $\sigma_{ii}$ and the conditional variances $\omega_{jj}$. It will
be shown that the convergence rate of CLIME obtained by Cai, Liu and Luo
(2011) is not optimal. The quantity $M_{n,p}$ is the upper bound of the
matrix $\ell_{1}$ norm which is unknown in practice. The cross validation
method can be used to choose the tuning parameter in CLIME. However, the
estimator obtained through CV can be variable and its theoretical properties
are unclear. In contrast, the ACLIME procedure proposed in the present paper
does not depend on any unknown parameters and it will be shown that the
estimator is minimax rate optimal.
\end{remark}

\section{Properties of ACLIME and Minimax Upper Bounds}
\label{optimality.sec}

We now study the properties of the ACLIME estimator $\hat{\Omega}
$ proposed in Section \ref{method.sec}. We shall begin with the Gaussian
case where $X\sim N(\mu ,\Sigma )$. Extensions to non-Gaussian distributions will be discussed later. The following result shows that the
ACLIME estimator adaptively attains the convergence rate of
\begin{equation*}
M_{n,p}^{1-q}c_{n,p}\left( \frac{\log p}{n}\right) ^{(1-q)/2}
\end{equation*}%
over the class of sparse precision matrices $\mathcal{G}%
_{q}(c_{n,p},M_{n,p}) $ defined in \eqref{parameterspace} under the matrix $%
\ell _{w}$ norm losses for all $1\leq w\leq \infty $. The lower bound given in
Section \ref{lowerbd.sec} shows that this rate is indeed optimal and thus
ACLIME adapts to both sparsity patterns and this class of loss functions.

\begin{theorem}
\label{cth1-1}
Suppose we observe a random sample $X_1, ..., X_n \stackrel{iid}{\sim} N_p(\mu, \Sigma)$. Let $\Omega=\Sigma^{-1}$ be the precision matrix. Let $\delta \geq 2$, $\log p=O(n^{1/3})$ and
\begin{equation}  \label{le4}
c_{n,p}=O\left( n^{\frac{1-q}{2}}/(\log p)^{\frac{1-q}{2}}\right).
\end{equation}%
Then for some constant $C>0$
\begin{equation*}
\inf_{\Omega\in \mathcal{G}_{q}(c_{n,p}, M_{n,p})}\mathbb{P}\left( \Vert
\hat{\Omega}-\Omega \Vert_w \leq CM_{n,p}^{1-q}c_{n,p}\left( \frac{\log p}{n}%
\right) ^{\frac{1-q}{2}}\right) \geq 1-O\left( (\log p)^{-{\frac{1}{2}} }p^{-%
{\frac{\delta^{2}}{4}}+1}\right)
\end{equation*}
for all $1\le w \le \infty$.
\end{theorem}

For $q=0$ a sufficient condition for estimating $\Omega $
consistently under the spectral norm is
\begin{equation*}
M_{n,p}c_{n,p}\sqrt{\frac{n}{\log p}}=o(1),  \quad{\rm i.e.,}\quad   M_{n,p}c_{n,p}=o\left( \sqrt{\frac{n}{\log p}}\right) .
\end{equation*}%
This implies that the total number of nonzero elements on each column needs
be $\ll \sqrt{n}$ in order for the precision matrix to be estimated
consistently over $\mathcal{G}_{0}(c_{n,p},M_{n,p})$. In Theorem \ref{Operlower.bd.thm} we show that the upper bound $M_{n,p}c_{n,p}\sqrt{\frac{\log p}{n}}$ is indeed rate optimal over $\mathcal{G}_{0}(c_{n,p},M_{n,p})$.

We now consider the rate of convergence under the expectation. For technical
reasons, we require the constant $\delta\geq 3$ in this case.


\begin{theorem}
\label{cminiupp}
Suppose we observe a random sample $X_1, ..., X_n \stackrel{iid}{\sim} N_p(\mu, \Sigma)$. Let $\Omega=\Sigma^{-1}$ be the precision matrix.
Let $\log p=O(n^{1/3})$ and $\delta \geq 3$. Suppose that $p\geq n^{13/(\delta^{2}-8)}$ and
\begin{equation*}
c_{n,q}=o((n/\log p)^{\frac{1}{2}-\frac{q}{2}}).
\end{equation*}%
The ACLIME estimator $\hat{\Omega}$ satisfies, for all $1\le w \le \infty$ and $0\le q < 1$,
\begin{equation*}
\sup_{\mathcal{G}_{q}(c_{n,p}, M_{n,p})}\mathbb{E}\left\Vert \hat{\Omega}%
-\Omega \right\Vert_w ^{2}\leq CM_{n,p}^{2-2q}c_{n,p}^{2}\left( \frac{\log p%
}{n}\right) ^{1-q},
\end{equation*}
for some constant $C>0$.
\end{theorem}

Theorem \ref{cminiupp} can be extended to non-Gaussian distributions.$\ $Let
$Z=\left( Z_{1},Z_{2},\ldots ,Z_{n}\right) ^{\prime }$ be a $p-$variate
random variable with mean $\mu $ and covariance matrix $\Sigma =\left(
\sigma _{ij}\right) _{1\leq i,j\leq p}$. Let $\Omega =\left( \omega _{ij}\right) _{1\leq i,j\leq p}$ be the precision matrix. Define $%
Y_{i}=(Z_{i}-\mu _{i})/\sigma _{ii}^{1/2}$, $1\leq i\leq p$ and $%
(W_{1},\ldots ,W_{p})^{\prime }:=\Omega (Z-\mu )$. Assume that there exist
some positive constants $\eta $ and $M$ such that for all $1\leq i\leq p$,
\begin{equation}
\mathbb{E}\exp (\eta Y_{i}^{2})\leq M,\quad \mathbb{E}\exp (\eta
W_{i}^{2}/\omega _{ii})\leq M\text{.\vspace{3mm}}  \label{T1}
\end{equation}
Then we have the following result.

\begin{theorem}
\label{Generalupp} Suppose we observe an i.i.d. sample $X_{1},...,X_{n}$
with the precision matrix $\Omega $ satisfying Condition (\ref{T1}). Let $\log
p=O(n^{1/3})$, $p\geq n^{\gamma}$ for some $\gamma>0$. Suppose that
\[
c_{n,q}=o((n/\log p)^{\frac{1}{2}-\frac{q}{2}}).
\]%
Then there is a $\delta$ depending only on $\eta$, $M$ and $\gamma$ such that the ACLIME estimator $\hat{\Omega}$ satisfies, for all $1\leq w\leq \infty $
and $0\leq q<1$,
\[
\sup_{\mathcal{G}_{q}(c_{n,p},M_{n,p})}\mathbb{E}\left\Vert \hat{\Omega}%
-\Omega \right\Vert _{w}^{2}\leq CM_{n,p}^{2-2q}c_{n,p}^{2}\left( \frac{\log
p}{n}\right) ^{1-q},
\]%
for some constant $C>0$.
\end{theorem}

\begin{remark}
Under Condition (\ref{T1}) it can be shown that an analogous result to
Lemma \ref{var.bnd.lem} in Section \ref{method.sec} holds with some $\delta$ depending only on $\eta$ and $M$. Thus, it can be proved that, under Condition (\ref{T1}),
Theorem \ref{Generalupp} holds. The proof is similar to
that of Theorem \ref{cminiupp}. A practical way to choose $\delta $ is
using cross validation.
\end{remark}
\begin{remark}
Theorems \ref{cth1-1}, \ref{cminiupp} and \ref{Generalupp} follow mainly from the convergence
rate under the element-wise $\ell _{\infty }$ norm and the inequality $\Vert
M\Vert _{w}\leq \Vert M\Vert _{1}$ for any symmetric matrix $M$ from Lemma \ref{lwl2}. The
convergence rate under element-wise norm plays an important role in
graphical model selection and in establishing the convergence rate under other
matrix norms, such as the Frobenius norm $\|\cdot\|_{F}$. Indeed, from the proof, Theorems \ref{cth1-1}, \ref{cminiupp} and \ref{Generalupp} hold under the matrix $\ell _{1}$
norm. More specifically, under the conditions of Theorems \ref{cminiupp} and \ref{Generalupp} we have
\begin{eqnarray*}
\sup_{\mathcal{G}_{q}(c_{n,p},M_{n,p})}\mathbb{E}|\hat{\Omega}-\Omega
|_{\infty }^{2} &\leq &CM_{n,p}^{2}\frac{\log p}{n}, \\
\sup_{\mathcal{G}_{q}(c_{n,p},M_{n,p})}\mathbb{E}\Vert \hat{\Omega}-\Omega
\Vert _{1}^{2} &\leq &CM_{n,p}^{2-2q}c_{n,p}^{2}\left( \frac{\log p}{n}%
\right) ^{1-q}, \\
\sup_{\mathcal{G}_{q}(c_{n,p},M_{n,p})}\frac{1}{p}\mathbb{E}\Vert \hat{\Omega%
}-\Omega \Vert _{F}^{2} &\leq &CM_{n,p}^{2-q}c_{n,p}\left( \frac{\log p}{n}%
\right) ^{1-q/2}.
\end{eqnarray*}
\end{remark}

\begin{remark}
\label{weak.lq.remark}
The results in this section can be easily extended to the weak $\ell _{q}$
ball with $0\leq q<1$ to model the sparsity of the precision matrix $\Omega $. A weak $\ell _{q}$ ball of radius $c$ in $R^{p}$ is defined as follows,
\begin{equation*}
B_{q}(c)=\left\{ \xi \in \mathbb{R}^{p}:\;|\xi |_{(k)}^{q}\leq ck^{-1},\quad %
\mbox{for all $k=1, ..., p$}\right\} ,
\end{equation*}%
where $|\xi |_{(1)}\geq $ $|\xi |_{(2)}\geq \ldots \geq $ $|\xi |_{(p)}$.
Let%
\beq
\label{weak.lq.ball}
\mathcal{G}_{q}^{\ast }(c_{n,p},M_{n,p})=\left\{
\begin{array}{c}
\Omega =\left( \omega _{ij}\right) _{1\leq i,j\leq p}:\omega _{\cdot ,j}\in
B_{q}(c_{n,p}), \\
\Vert \Omega \Vert _{1}\leq M_{n,p},\text{ }\lambda _{\max }(\Omega
)/\lambda _{\min }(\Omega )\leq M_{1},\Omega \succ 0%
\end{array}%
\right\} .
\eeq
Theorems \ref{cth1-1}, \ref{cminiupp} and \ref{Generalupp} hold with the parameter space $%
G_{q}(c_{n,p},M_{n,p})$ replaced by $G_{q}^{\ast }(c_{n,p},M_{n,p})$ by a
slight extension of Lemma \ref{le5} for the $\ell _{q}$ ball to for the weak
$\ell _{q}$ ball similar to Equation (51) in Cai and Zhou (2012).
\end{remark}

\section{Minimax Lower Bounds}

\label{lowerbd.sec} 

Theorem \ref{cminiupp} shows that the ACLIME estimator adaptively attains
the rate of convergence
\begin{equation}
M_{n,p}^{2-2q}c_{n,p}^{2}\left( \frac{\log p}{n}\right) ^{1-q}
\label{opt.rate}
\end{equation}%
under the squared matrix $\ell _{w}$ norm loss for $1\leq w\leq \infty $
over the collection of the parameter spaces $\mathcal{G}_{q}(c_{n,p},M_{n,p})$.
In this section we shall show that the rate of convergence given in %
\eqref{opt.rate} cannot be improved by any other estimator and thus is
indeed optimal among all estimators by establishing minimax lower bounds for
estimating sparse precision matrices under the squared matrix $\ell
_{w}$ norm.


\begin{theorem}
\label{Operlower.bd.thm} Let $X_{1},\ldots ,X_{n}\overset{iid}{\sim }N_p(\mu,\Sigma)$ with $p>c_{1}n^{\beta }$ for some constants $\beta >1$ and $c_{1}>0$.  Assume that
\begin{equation}
cM_{n,p}^{q}\left( {\log p\over n}\right) ^{q\over 2}\leq c_{n,p}=o\left( M_{n,p}^{q}n^{\frac{1-q}{2}}\left( \log p\right) ^{-\frac{3-q}{2}}\right)
\label{cnp}
\end{equation}%
for some constant $c>0$. The minimax risk for estimating the precision matrix $\Omega =\Sigma^{-1}$ over the parameter space $\mathcal{G}_{q}(c_{n,p},M_{n,p})$ under
the condition (\ref{cnp}) satisfies
\begin{equation*}
\inf_{\hat{\Omega}}\sup_{\mathcal{G}_{q}(c_{n,p},M_{n,p})}\mathbb{E}%
\left\Vert \hat{\Omega}-\Omega \right\Vert _{w}^{2}\geq CM_{n,p}^{2-2q}c
_{n,p}^{2}\left( \frac{\log p}{n}\right) ^{1-q}
\end{equation*}%
for some constant $C>0$ and for all $1\leq w\leq \infty $.
\end{theorem}

The proof of Theorem \ref{Operlower.bd.thm} is involved. We shall discuss
the key technical tools and outline the important steps in the proof of
Theorem \ref{Operlower.bd.thm} in this section. The detailed proof is given
in Section \ref{proof.sec}.


\subsection{A General Technical Tool}

\label{Glowbd.sec}


We use a lower bound technique introduced in Cai and Zhou (2012), which is
particularly well suited for treating \textquotedblleft
two-directional\textquotedblright\ problems such as matrix estimation. The
technique can be viewed as a generalization of both Le Cam's method and
Assouad's Lemma, two classical lower bound arguments. Let $X$ be an
observation from a distribution $\mathbb{P}_{\theta }$ where $\theta $
belongs to a parameter set $\Theta $ which has a special tensor structure.
For a given positive integer $r$ and a finite set $B\subset \mathbb{R}%
^{p}/\left\{ 0_{1\times p}\right\} $, let $\Gamma =\left\{ 0,1\right\} ^{r}$
and $\Lambda \subseteq B^{r}$. Define
\begin{equation}
\Theta =\Gamma \otimes \Lambda =\left\{ \left( \gamma ,\lambda \right)
:\gamma \in \Gamma \text{ and }\lambda \in \Lambda \right\} \text{.}
\label{theta}
\end{equation}%
In comparison, the standard lower bound arguments work with either $\Gamma $
or $\Lambda $ alone. For example, the Assouad's Lemma considers only the
parameter set $\Gamma $ and the Le Cam's method typically applies to a
parameter set like $\Lambda $ with $r=1$. Cai and Zhou (2012) gives a lower
bound for the maximum risk over the parameter set $\Theta $ to the problem
of estimating a functional $\psi (\theta )$, belonging to a metric space
with metric $d$.

We need to introduce a few notations before formally stating the lower
bound. For two distributions $\mathbb{P}$ and $\mathbb{Q}$ with densities $p$
and $q$ with respect to any common dominating measure $\mu $, the total
variation affinity is given by $\Vert \mathbb{P}\wedge \mathbb{Q}\Vert =\int
p\wedge qd\mu $. For a parameter $\gamma =(\gamma _{1},...,\gamma _{r})\in
\Gamma $ where $\gamma _{i}\in \{0,1\}$, define
\begin{equation}
H\left( \gamma ,\gamma ^{\prime }\right) =\sum_{i=1}^{r}\left\vert \gamma
_{i}-\gamma _{i}^{\prime }\right\vert   \label{H}
\end{equation}%
be the Hamming distance on $\left\{ 0,1\right\} ^{r}$.

Let $D_{\Lambda }=$\textrm{Card}$\left( \Lambda \right) $. For a given $a\in
\{0,1\}$ and $1\leq i\leq r$, we define the mixture distribution $\mathbb{%
\bar{P}}_{a,i}$ by
\begin{equation}
\mathbb{\bar{P}}_{a,i}=\frac{1}{2^{r-1}D_{\Lambda }}\sum_{\theta }\{\mathbb{P%
}_{\theta }:\;\gamma _{i}(\theta )=a\}.  \label{avepi}
\end{equation}%
So $\mathbb{\bar{P}}_{a,i}$ is the mixture distribution over all $P_{\theta }
$ with $\gamma _{i}(\theta )$ fixed to be $a$ while all other components of $%
\theta $ vary over all possible values. In our construction of the parameter
set for establishing the minimax lower bound, $r$ is the number of possibly
non-zero rows in the upper triangle of the covariance matrix, and $\Lambda $
is the set of matrices with $r$ rows to determine the upper triangle matrix.

\begin{lemma}
\label{AL} For any estimator $T$ of $\psi (\theta )$ based on an observation
from the experiment $\left\{ \mathbb{P}_{\theta },\theta \in \Theta \right\}
$, and any $s>0$%
\begin{equation}
\max_{\Theta }2^{s}\mathbb{E}_{\theta }d^{s}\left( T,\psi \left( \theta
\right) \right) \geq \alpha \frac{r}{2}\min_{1\leq i\leq r}\left\Vert
\mathbb{\bar{P}}_{0,i}\wedge \mathbb{\bar{P}}_{1,i}\right\Vert
\label{rhsLemma1}
\end{equation}%
where $\mathbb{\bar{P}}_{a,i}$ is defined in Equation (\ref{avepi}) and $%
\alpha $ is given by
\begin{equation}
\alpha =\min_{\left\{ (\theta ,\theta ^{\prime }):H(\gamma (\theta ),\gamma
(\theta ^{\prime }))\geq 1\right\} }\frac{d^{s}(\psi (\theta ),\psi (\theta
^{\prime }))}{H(\gamma (\theta ),\gamma (\theta ^{\prime }))}\text{.}
\label{alpha}
\end{equation}
\end{lemma}

We introduce some new notations to study the affinity $\left\Vert
\mathbb{\bar{P}}_{0,i}\wedge \mathbb{\bar{P}}_{1,i}\right\Vert $ in Equation
(\ref{rhsLemma1}). Denote the projection of $\theta \in \Theta $ to $\Gamma $
by $\gamma \left( \theta \right) =\left( \gamma _{i}\left( \theta \right)
\right) _{1\leq i\leq r}$ and to $\Lambda $ by $\lambda \left( \theta
\right) =\left( \lambda _{i}\left( \theta \right) \right) _{1\leq i\leq r}$.
More generally we define $\gamma _{A}\left( \theta \right) =\left( \gamma
_{i}\left( \theta \right) \right) _{i\in A}$ for a subset $A\subseteq
\left\{ 1,2,\ldots ,r\right\} $, a projection of $\theta $ to a subset of $%
\Gamma $. A particularly useful example of set $A$ is
\begin{equation*}
\left\{ -i\right\} =\left\{ 1,\ldots ,i-1,i+1,\cdots ,r\right\} ,
\end{equation*}%
for which $\gamma _{-i}\left( \theta \right) =\left( \gamma _{1}\left(
\theta \right) ,\ldots ,\gamma _{i-1}\left( \theta \right) ,\gamma
_{i+1}\left( \theta \right) ,\gamma _{r}\left( \theta \right) \right) $. $%
\lambda _{A}\left( \theta \right) $ and $\lambda _{-i}\left( \theta \right) $
are defined similarly. We denote the set $\left\{ \lambda _{A}\left( \theta
\right) :\theta \in \Theta \right\} $ by $\Lambda _{A}$. For $a\in \left\{
0,1\right\} $, $b\in \left\{ 0,1\right\} ^{r-1}$, and $c\in \Lambda
_{-i}\subseteq B^{r-1}$, let%
\begin{equation*}
D_{\Lambda _{i}\left( a,b,c\right) }=\mathrm{Card}\left\{ \gamma \in \Lambda
:\gamma _{i}(\theta )=a,\gamma _{-i}(\theta )=b\text{ and }\lambda
_{-i}(\theta )=c\right\}
\end{equation*}%
and define
\begin{equation}
\mathbb{\bar{P}}_{\left( a,i,b,c\right) }=\frac{1}{D_{\Lambda _{i}\left(
b,c\right) }}\sum_{\theta }\{\mathbb{P}_{\theta }:\;\gamma _{i}(\theta
)=a,\gamma _{-i}(\theta )=b\text{ and }\lambda _{-i}(\theta )=c\}\text{.}
\label{avepibd}
\end{equation}%
In other words, $\mathbb{\bar{P}}_{(a,i,b,c)}$ is the mixture distribution
over all $\mathbb{P}_{\theta }$ with $\lambda _{i}(\theta )$ varying over
all possible values while all other components of $\theta $ remain fixed.

The following lemma gives a lower bound for the affinity in Equation
(\ref{rhsLemma1}). See SEction 2 of  Cai and Zhou (2012) for more details.

\begin{lemma}
\label{aff}Let $\mathbb{\bar{P}}_{a,i}$ and $\mathbb{\bar{P}}_{\left(
a,i,b,c\right) }$ be defined in Equation (\ref{avepi})and (\ref{avepibd})
respectively, then
\begin{equation*}
\left\Vert \mathbb{\bar{P}}_{0,i}\wedge \mathbb{\bar{P}}_{1,i}\right\Vert
\geq \underset{\gamma _{-i},\lambda _{-i}}{\mathrm{Average}}\left\Vert
\mathbb{\bar{P}}_{\left( 0,i,\gamma _{-i},\lambda _{-i}\right) }\wedge
\left( \mathbb{\bar{P}}_{\left( 1,i,\gamma _{-i},\lambda _{-i}\right)
}\right) \right\Vert \text{,}
\end{equation*}
where the average over $\gamma _{-i}$ and $\lambda _{-i}$ is induced by the uniform distribution
over $\Theta$.
\end{lemma}


\subsection{Lower Bound for Estimating Sparse Precision Matrix}


We now apply the lower bound technique developed in Section \ref{Glowbd.sec}
to establish rate sharp results under the matrix $\ell _{w}$ norm. Let $%
X_{1},\ldots ,X_{n}\overset{iid}{\sim }N_p(\mu ,\Omega^{-1})$
with $p>c_{1}n^{\beta }$ for some $\beta >1$ and $c_{1}>0$, where $\Omega
\in $ $\mathcal{G}_{q}(c_{n,p},M_{n,p})$. The proof of Theorem \ref{Operlower.bd.thm} contains four
major steps. We first reduce the minimax lower bound under the general
matrix $\ell _{w}$ norm, $1\leq w\leq \infty $, to under the spectral norm.
In the second step we construct in detail a subset $\mathcal{F}_{\ast }$ of
the parameter space $\mathcal{G}_{q}(c_{n,p},M_{n,p})$ such that the
difficulty of estimation over $\mathcal{F}_{\ast }$ is essentially the same
as that of estimation over $\mathcal{G}_{q}(c_{n,p},M_{n,p})$, the third
step is the application of Lemma \ref{AL} to the carefully constructed
parameter set, and finally in the fourth step we calculate the factors $%
\alpha $ defined in (\ref{alpha}) and the total variation affinity between
two multivariate normal mixtures. We outline the main ideas of the proof
here and leave detailed proof of some technical results to Section \ref%
{proof.sec}.

\medskip\noindent \textbf{Proof of Theorem \ref{Operlower.bd.thm}:} We shall
divide the proof into four major steps.

\medskip\noindent \textbf{Step 1: Reducing the general problem to the lower
bound under the spectral norm. }The following lemma implies that the minimax
lower bound under the spectral norm yields a lower bound under the general
matrix $\ell _{w}$ norm up to a constant factor $4$.

\begin{lemma}
\label{reduction}Let $X_{1},\ldots ,X_{n}\overset{iid}{\sim }N(\mu ,\Omega^{-1})$, and $\mathcal{F}$ be any parameter space of precision
matrices. The minimax risk for estimating the precision matrix $%
\Omega $ over $\mathcal{F}$ satisfies
\begin{equation}
\inf_{\hat{\Omega}}\sup_{\mathcal{F}}\mathbb{E}\left\Vert \hat{\Omega}%
-\Omega \right\Vert _{w}^{2}\geq \frac{1}{4}\inf_{\hat{\Omega}}\sup_{%
\mathcal{F}}\mathbb{E}\left\Vert \hat{\Omega}-\Omega \right\Vert _{2}^{2}
\label{symmest}
\end{equation}%
for all $1\leq w\leq \infty $.
\end{lemma}

\medskip \noindent \textbf{Step 2: Constructing the parameter set.} Let $%
r=\lceil p/2\rceil $ and let $B$ be the collection of all vectors $\left(
b_{j}\right) _{1\leq j\leq p}$ such that $b_{j}=0$ for $1\leq j\leq p-r$ and
$b_{j}=0$ or $1$ for $p-r+1\leq j\leq p$ under the constraint $\left\Vert
b\right\Vert _{0}=k$ (to be defined later). For each $b\in B$ and each $%
1\leq m\leq r$, define a $p\times p$ matrix $\lambda _{m}(b)$ by making the $%
m$th row of $\lambda _{m}(b)$ equal to $b$ and the rest of the entries $0$.
It is clear that \textrm{Card}($B$)$=\binom{r}{k}$. Set $\Gamma =\left\{
0,1\right\} ^{r}$. Note that each component $b_{i}$ of $\lambda
=(b_{1},...,b_{r})\in \Lambda $ can be uniquely associated with a $p\times p$
matrix $\lambda _{i}(b_{i})$. $\Lambda $ is the set of all matrices $\lambda
$ with the every column sum less than or equal to $2k$. Define $\Theta
=\Gamma \otimes \Lambda $ and let $\epsilon _{n,p}\in \mathbb{R}$ be fixed.
(The exact value of $\epsilon _{n,p}$ will be chosen later.) For each $%
\theta =(\gamma ,\lambda )\in \Theta $ with $\gamma =(\gamma _{1},...,\gamma
_{r})$ and $\lambda =(b_{1},...,b_{r})$, we associate $\theta $ with a
precision matrix $\Omega (\theta )$ by
\begin{equation*}
\Omega (\theta )=\frac{M_{n,p}}{2}\left[ I_{p}+\epsilon
_{n,p}\sum_{m=1}^{r}\gamma _{m}\lambda _{m}(b_{m})\right] .
\end{equation*}%
Finally we define a collection $\mathcal{F}_{\ast }$ of precision
matrices as
\begin{equation*}
\mathcal{F}_{\ast }=\left\{ \Omega (\theta ):\Omega (\theta )=\frac{M_{n,p}}{%
2}\left[ I_{p}+\epsilon _{n,p}\sum_{m=1}^{r}\gamma _{m}\lambda _{m}(b_{m})%
\right] ,\;\theta =(\gamma ,\lambda )\in \Theta \right\} .
\end{equation*}%
We now specify the values of $\epsilon _{n,p}$ and $k$. Set%
\begin{equation}
\epsilon _{n,p}=\upsilon \sqrt{\frac{\log p}{n}}\text{, for some }0<\upsilon
<\min \left\{ \left( \frac{c}{2}\right) ^{1/q}\text{, }\frac{\beta -1}{%
8\beta }\right\} ,  \label{epsilonnp}
\end{equation}%
and
\begin{equation}
{k=\left\lceil 2^{-1}c_{n,p}(M_{n,p}\epsilon _{n,p})^{-q}\right\rceil -1}%
\text{{.}}  \label{kdef}
\end{equation}%
which is at least $1$ from Equation (\ref{epsilonnp}). Now we show $\mathcal{%
F}_{\ast }$ is a subset of the parameter space $\mathcal{G}%
_{q}(c_{n,p},M_{n,p})$. From the definition of $k$ in (\ref{kdef}) note that
\begin{equation}
\max_{j\leq p}\sum_{i\neq j}|\omega _{ij}|^{q}\leq 2\cdot 2^{-1}\rho
_{n,p}\left( M_{n,p}\epsilon _{n,p}\right) ^{-q}\cdot \left( \frac{M_{n,p}}{2%
}\epsilon _{n,p}\right) ^{q}\leq c_{n,p}.  \label{qnorm}
\end{equation}%
From Equation (\ref{cnp}) we have $c_{n,p}=o\left( M_{n,p}^{q}n^{\frac{1-q}{2%
}}\left( \log p\right) ^{-\frac{3-q}{2}}\right) $, which implies
\begin{equation}
2k\epsilon _{n,p}\leq c_{n,p}\epsilon _{n,p}^{1-q}M_{n,p}^{-q}=o\left(
1/\log p\right) ,  \label{kepsilon}
\end{equation}%
then%
\begin{equation}
\max_{i}\sum_{j}|\omega _{ij}|\leq \frac{M_{n,p}}{2}\left( 1+2k\epsilon
_{n,p}\right) \leq M_{n,p}.  \label{1norm}
\end{equation}%
Since $\left\Vert A\right\Vert _{2}\leq \left\Vert A\right\Vert _{1}$, we
have%
\begin{equation*}
\left\Vert \epsilon _{n,p}\sum_{m=1}^{r}\gamma _{m}\lambda
_{m}(b_{m})\right\Vert _{2}\leq \left\Vert \epsilon
_{n,p}\sum_{m=1}^{r}\gamma _{m}\lambda _{m}(b_{m})\right\Vert _{1}\leq
2k\epsilon _{n,p}=o\left( 1\right) ,
\end{equation*}%
which implies that every $\Omega (\theta )$ is diagonally dominant and
positive definite, and
\begin{equation}
\lambda _{\max }\left( \Omega \right) \leq \frac{M_{n,p}}{2}\left(
1+2k\epsilon _{n,p}\right) \text{, and }\lambda _{\min }\left( \Omega
\right) \geq \frac{M_{n,p}}{2}\left( 1-2k\epsilon _{n,p}\right)
\label{posi}
\end{equation}%
which immediately implies%
\begin{equation}
\frac{\lambda _{\max }\left( \Omega \right) }{\lambda _{\min }\left( \Omega
\right) }<M_{1}\text{.}  \label{cond}
\end{equation}%
Equations (\ref{qnorm}), (\ref{1norm}), (\ref{posi}) and (\ref{cond}) all
together imply $\mathcal{F}_{\ast }\subset \mathcal{G}_{q}(c_{n,p},M_{n,p})$.

\medskip\noindent \textbf{Step 3: Applying the general lower bound argument.}
Let $X_{1},\ldots ,X_{n}\overset{iid}{\sim }N_p\left( 0,\left( \Omega (\theta
)\right) ^{-1}\right) $ with $\theta \in \Theta $ and denote the joint
distribution by $P_{\theta }$. Applying Lemmas \ref{AL} and \ref{aff} to the
parameter space $\Theta $, we have
\begin{equation}
\inf_{\hat{\Omega}}\max_{\theta \in \Theta }2^{2}E_{\theta }\left\Vert \hat{%
\Omega}-\Omega (\theta )\right\Vert _{2}^{2}\geq \alpha \cdot \frac{p}{4}%
\cdot \min_{i}\underset{\gamma _{-i},\lambda _{-i}}{\mathrm{Average}}%
\left\Vert \mathbb{\bar{P}}_{\left( 0,i,\gamma _{-i},\lambda _{-i}\right)
}\wedge \mathbb{\bar{P}}_{\left( 1,i,\gamma _{-i},\lambda _{-i}\right)
}\right\Vert  \label{lowerbound*}
\end{equation}%
where
\begin{equation}
\alpha =\min_{\left\{ (\theta ,\theta ^{\prime }):H(\gamma (\theta ),\gamma
(\theta ^{\prime }))\geq 1\right\} }\frac{\left\Vert \Omega (\theta )-\Omega
(\theta ^{\prime })\right\Vert _{2}^{2}}{H(\gamma (\theta ),\gamma (\theta
^{\prime }))}  \label{alpha1}
\end{equation}%
and $\mathbb{\bar{P}}_{0,i}$ and $\mathbb{\bar{P}}_{1,i}$ are defined as in (%
\ref{avepi}).

\medskip\noindent \textbf{Step 4: Bounding the per comparison loss $\alpha $
defined in (\ref{alpha1}) and the affinity $\underset{i}{\min}\underset{%
\gamma _{-i},\lambda _{-i}}{\mathrm{Average}}\left\Vert \mathbb{\bar{P}}%
_{\left( 0,i,\gamma _{-i},\lambda _{-i}\right) }\wedge \mathbb{\bar{P}}%
_{\left( 1,i,\gamma _{-i},\lambda _{-i}\right) }\right\Vert$ in (\ref%
{lowerbound*}).} This is done separately in the next two lemmas which are
proved in detailed in Section \ref{proof.sec}.

\begin{lemma}
\label{dffbd} The per comparison loss $\alpha $ defined in (\ref{alpha1})
satisfies
\begin{equation*}
\alpha \geq \frac{(M_{n,p}k\epsilon _{n,p})^{2}}{4p}.
\end{equation*}
\end{lemma}

\begin{lemma}
\label{affbd} Let $X_{1},\ldots ,X_{n}\overset{iid}{\sim }N\left( 0,\left(
\Omega (\theta )\right) ^{-1}\right) $ with $\theta \in \Theta $ and denote
the joint distribution by $\mathbb{P}_{\theta }$. For $a\in \{0,1\}$ and $%
1\leq i\leq r$, define $\mathbb{\bar{P}}_{\left( a,i,b,c\right) }$ as in (%
\ref{avepibd}). Then there exists a constant $c_{1}>0$ such that
\begin{equation*}
\min_{i}\underset{\gamma _{-i},\lambda _{-i}}{\mathrm{Average}}\left\Vert
\mathbb{\bar{P}}_{\left( 0,i,\gamma _{-i},\lambda _{-i}\right) }\wedge
\mathbb{\bar{P}}_{\left( 1,i,\gamma _{-i},\lambda _{-i}\right) }\right\Vert
\geq c_{1}.
\end{equation*}
\end{lemma}

Finally, the minimax lower bound for estimating a sparse precision matrix
over the collection $\mathcal{G}_{q}(c_{n,p},M_{n,p})$ is obtained by
putting together (\ref{lowerbound*}) and Lemmas \ref{dffbd} and \ref{affbd},%
\begin{eqnarray*}
\inf_{\hat{\Omega}}\sup_{\mathcal{G}_{q}(c_{n,p},M_{n,p})}\mathbb{E}%
\left\Vert \hat{\Omega}-\Omega (\theta )\right\Vert _{2}^{2} &\geq
&\max_{\Omega (\theta )\in \mathcal{F}_{\ast }}E_{\theta }\left\Vert \hat{%
\Omega}-\Omega (\theta )\right\Vert _{2}^{2}\geq \frac{\left(
M_{n,p}k\epsilon _{n,p}\right) ^{2}}{4p}\cdot \frac{p}{16}\cdot c_{1} \\
&\geq &\frac{c_{1}}{64}(M_{n,p}k\epsilon _{n,p})^{2}=c_{2}{M_{n,p}^{2-2q}}%
c_{n,p}^{2}\left( \frac{\log p}{n}\right) ^{1-q}\text{,}
\end{eqnarray*}%
for some constant $c_{2}>0$. \qed

Putting together the minimax upper and lower bounds in Theorems \ref{cminiupp} and \ref{Operlower.bd.thm} as well as Remark \ref{weak.lq.remark} yields the optimal rates of convergence for estimating $\Omega$ over the collection of the $\ell_q$ balls $\mathcal{G}_{q}(c_{n,p},M_{n,p})$ defined in (\ref{parameterspace}) as well as the collection of the weak $\ell_q$ balls $\mathcal{G}_{q}^*(c_{n,p},M_{n,p})$ defined in \eqref{weak.lq.ball}.

\begin{theorem}
\label{MinimaxOpe} Suppose we observe a random sample $X_{i}{\overset{\text{%
i.i.d.}}{\sim }}N_{p}(\mu ,\Sigma )$, $i=1,2,\ldots ,n$. Let $\Omega =\Sigma
^{-1}$ be the precision matrix. Assume that $\log p=O(n^{1/3})$ and
\begin{equation}
cM_{n,p}^{q}\left( \frac{\log p}{n}\right) ^{q\over 2} \leq c_{n,p}=o\left( M_{n,p}^{q}n^{%
\frac{1-q}{2}}\left( \log p\right) ^{-\frac{3-q}{2}}\right)   \label{cond.c}
\end{equation}%
for some constant $c>0$. Then
\begin{equation}
\inf_{\hat{\Omega}}\sup_{\Omega\in {\cal G}}\mathbb{E}%
\left\Vert \hat{\Omega}-\Omega \right\Vert _{w}^{2}\asymp
M_{n,p}^{2-2q}c_{n,p}^{2}\left( \frac{\log p}{n}\right) ^{1-q}
\label{rateOper}
\end{equation}%
for all $1\leq w\leq \infty $, where ${\cal G} = \mathcal{G}_{q}(c_{n,p},M_{n,p})$ or $\mathcal{G}_{q}^*(c_{n,p},M_{n,p})$.
\end{theorem}

\section{Numerical results}
\label{simulation.sec}

In this section, we consider the numerical performance of ACLIME. In particular, we shall  compare  the performance of ACLIME with that of  CLIME.  The following three graphical models are considered. Let $D={\rm diag}(U_{1},\ldots,U_{p})$, where $U_{i}$, $1\leq i\leq p$, are i.i.d. uniform random variables on the interval $(1,5)$. Let $\Sigma=\Omega^{-1}=D\Omega_{1}^{-1} D$. The matrix $D$ makes the diagonal entries in $\Sigma$ and $\Omega$ different.

\begin{itemize}

\item{ {\bf Band graph}. Let $\Omega_{1}=(\omega_{ij})$, where  $\omega_{ii}=1$, $\omega_{i,i+1}=\omega_{i+1,i}=0.6$, $\omega_{i,i+2}=\omega_{i+2,i}=0.3$, $\omega_{ij}=0$ for $|i-j|\geq 3$. }

\item{{\bf AR(1) model}. Let $\Omega_{1}=(\omega_{ij})$, where $\omega_{ij}=(0.6)^{|j-i|}$.}

\item{{\bf Erd\"{o}s-R\'{e}nyi random graph}. Let $\Omega_{2}=(\omega_{ij})$, where $\omega_{ij}=u_{ij}*\delta_{ij}$,  $\delta_{ij}$ is the
Bernoulli random variable with success probability 0.05 and $u_{ij}$ is uniform random variable with distribution $U(0.4,0.8)$. We let $\Omega_{1}=\Omega_{2}+(|\min(\lambda_{\min})|+0.05)I_{p}$. It is easy to check that the matrix $\Omega_1$ is symmetric and  positive definite.}

\end{itemize}

We generate $n=200$ random  training samples from $N_p(0,\Sigma)$ distribution for $p=50,100,200$.  For ACLIME, we set $\delta=2$ in Step 1 and choose $\delta$ in
Step 2 by a cross validation method. To this end, we  generate an additional 200 testing samples. The tuning parameter in CLIME is selected by cross validation. Note that ACLIME chooses different tuning parameters for different columns and CLIME chooses a universal tuning parameter. The log-likehood loss
\begin{eqnarray*}
L(\hat{\Sigma}_{1},\Omega)=\log(\det(\Omega))-\langle\hat{\Sigma}_{1},\Omega\rangle,
\end{eqnarray*}
where $\hat{\Sigma}_{1}$ is the sample covariance matrix of the testing samples, is used in the cross validation method. For $\delta$ in (\ref{lambdan}), we let $\delta=\delta_{j}=j/50$, $1\leq j\leq 100$. For each $\delta_{j}$,
 ACLIME $\hat{\Omega}(\delta_{j})$ is obtained and the tuning parameter $\delta$ in (\ref{lambdan}) is selected by minimizing the following log-likehood loss
\begin{eqnarray*}
\hat{\delta}=\hat{j}/50,\mbox{~~where~~}\hat{j}=\argmin_{1\leq j\leq 100} L(\hat{\Sigma}_{1},\hat{\Omega}(\delta_{j})).
\end{eqnarray*}
The tuning parameter $\lambda_{n}$ in CLIME is also selected by cross validation. The detailed steps can be found in Cai, Liu and Luo (2011).

The empirical errors of ACLIME and CLIME estimators under various settings are summarized  in Table \ref{tb:simu} below. Three losses under the spectral norm, matrix $\ell_{1}$ norm and  Frobenius norm are given to compare the performance between ACLIME and CLIME. As can be seen from Table \ref{tb:simu},  ACLIME, which is tuning-free, outperforms CLIME in most of the cases for each of the three graphs.

\begin{table}[hptb]\addtolength{\tabcolsep}{1pt}
  \begin{center}
    \renewcommand{\arraystretch}{1.2}
    \begin{tabular}{|c|ccc|ccc|}
      \hline
      & \multicolumn{3}{c }{ACLIME}&\multicolumn{3}{c|}{CLIME} \\
      \hline
    $p$  &  50 &100& 200& 50&  100 & 200\\
      \hline
      & \multicolumn{6}{c|}{Spectral norm} \\
      \hline
      Band     &   0.30(0.01)&   0.45(0.01)  &  0.65(0.01)                               &0.32(0.01)    &0.50(0.01)     &    0.72(0.01)   \\
      \hline
     AR(1)     &   0.75(0.01)&   1.04(0.01)  &  1.25(0.01)                               &0.73(0.01)    &1.05(0.01)     &     1.30(0.01)   \\
      \hline
      E-R    &   0.65(0.03)&   0.95(0.02)  &  2.62(0.02)                               &0.72(0.03)    &1.21(0.04)     &     2.28(0.02)   \\
      \hline
      & \multicolumn{6}{c|}{Matrix $\ell_{1}$ norm} \\
      \hline
       Band     &   0.62(0.02)&   0.79(0.01)  &  0.94(0.01)                              &0.65(0.02)    &0.86(0.02)    &     0.99(0.01)   \\
             \hline
             AR(1)     &   1.19(0.02)&   1.62(0.02)  &  1.93(0.01)                               &1.17(0.01)    &1.59(.01)     &     1.89(0.01)   \\
      \hline
      E-R    &   1.47(0.08)&   2.15(0.06)  &  5.47(0.05)                               &1.53(0.06)    &2.34(0.06)     &     5.20(0.04)   \\
      \hline
      &\multicolumn{6}{c|}{Frobenius norm} \\
      \hline
       Band     &   0.80(0.01)&   1.61(0.02)  &  3.11(0.02)                               &0.83(0.01)    &1.73(0.02)     &     3.29(0.03)   \\
              \hline
              AR(1)     &   1.47(0.02)&   2.73(0.01)  &  4.72(0.01)                               &1.47(0.02)    &2.82(0.02)     &     4.97(0.01)   \\
      \hline
      E-R    &   1.53(0.05)&   3.15(0.03)  &  9.89(0.07)                              &1.62(0.04)    &3.61(0.05)     &    8.86(0.04)   \\
       \hline
    \end{tabular}
        \caption{Comparisons of ACLIME and CLIME for the three graphical models under three matrix norm losses. Inside the parentheses are the standard deviations of the empirical errors over 100 replications.}
    \label{tb:simu}
  \end{center}
\end{table}

\section{Discussions}
\label{discussion.sec}

We established in this paper  the optimal rates of convergence and introduced an adaptive method for  estimating sparse precision matrices under the matrix $\ell_w$ norm losses for $1\le w \le \infty$. The minimax rate of convergence under the Frobenius norm loss can also be easily established. As seen in the proof of Theorems \ref{cth1-1} and  \ref{cminiupp}, with probability tending to one,
 \begin{equation}
|\hat{\Omega}-\Omega |_{\infty }\leq CM_{n,p}\sqrt{\frac{\log p}{n}},
\label{infinitybound}
\end{equation}
for some constant $C>0$. From Equation (\ref{infinitybound}) one can immediately obtain the following risk upper bound under the Frobenius norm, which can be shown to  be rate optimal using a similar proof to that of Theorem \ref{Operlower.bd.thm}.
\begin{theorem}
\label{MinimaxF} Suppose we observe a random sample $X_{i}{\overset{\text{%
i.i.d.}}{\sim }}N_{p}(\mu ,\Sigma )$, $i=1,2,\ldots ,n$. Let $\Omega =\Sigma
^{-1}$ be the precision matrix. Under the assumption (\ref{cond.c}), the
minimax risk of estimating the precision matrix $\Omega $ over the class $%
\mathcal{G}_{q}(c_{n,p},M_{n,p})$ defined in (\ref{parameterspace}) satisfies%
\[
\inf_{\hat{\Omega}}\sup_{\mathcal{G}_{q}(c_{n,p},M_{n,p})}\mathbb{E}%
\frac{1}{p}\left\Vert \hat{\Omega}-\Omega \right\Vert _{F}^{2}\asymp
M_{n,p}^{2-q}c_{n,p}\left( \frac{\log p}{n}\right) ^{1-{\frac{q}{2}}}\text{.}
\]
\end{theorem}

As shown in Theorem \ref{MinimaxOpe}, the optimal rate of convergence for  estimating sparse precision matrices under the squared $\ell_w$ norm loss is $M_{n,p}^{2-2q}c_{n,p}^2\left( \frac{\log p}{n}\right) ^{1-q}$. It is interesting to compare this with the minimax rate of convergence for  estimating sparse covariance matrices under the same loss which is $c_{n,p}^2\left( \frac{\log p}{n}\right) ^{1-q}$ (cf. Theorem 1 in Cai and Zhou (2012)). These two convergence rates are similar, but have an important distinction. The difficulty of  estimating a sparse covariance matrix does not depend on the $\ell_1$ norm bound $M_{n,p}$, while the difficulty of estimating a sparse precision matrix does.

 As mentioned in the introduction, an important related problem to the estimation of precision matrix is  the recovery of a Gaussian graph which is equivalent to the estimation of the support of $\Omega$. Let $G=(V,E)$ be an undirected graph representing the conditional independence relations between the components of a random vector $X$. The vertex set $V$ contains the components of $X$, $V=X=\{V_{1},\dotsc ,V_{p}\}$. The edge set $E$ consists of ordered pairs $(i,j)$, indicating conditional dependence between the components $V_{i}$ and $V_{j}$. An edge between $V_{i}$ and $V_{j}$ is in the set $E$, i.e., $(i,j)\in E$, if and only  $\omega_{ij}=0$. The adaptive CLIME estimator, with an additional thresholding
step, can recover the support of $\Omega $. Define the estimator of the support of $\Omega$ by
\begin{eqnarray*}
\widehat{\text{SUPP}(\Omega)}=\{(i,j): |\hat{\omega}_{ij}|\geq \tau_{ij}\},
\end{eqnarray*}
where the choice of $\tau_{ij}$ depends on the bound $|\hat{\omega}_{ij}-\omega_{ij}|$. Equation (\ref{infinitybound}) implies that the right threshold levels are $\tau_{ij}=CM_{n,p}\sqrt{\log p/n}$. If the
magnitudes of the nonzero entries exceed $2CM_{n,p}\sqrt{\log p/n}$, then $\widehat{\text{SUPP}(\Omega)}$ recovers the support of $\Omega$ exactly.
In the context of covariance matrix estimation, Cai and Liu (2011) introduced an adaptive entry-dependent thresholding procedure to recover the support of $\Sigma$. That method is based on the sharp bound
\[
\max_{1\leq i\leq j\leq p}|\hat{\sigma}_{ij}-\sigma_{ij}|\leq 2\sqrt{\hat{\theta}_{ij}\log p/n},
\]
where $\hat{\theta}_{ij}$ is an estimator of Var$((X_{i}-\mu_{i})(X_{j}-\mu_{j}))$. It is natural to ask whether one can use data and entry-dependent threshold levels $\tau_{ij}$  to recover the support of $\Omega$. It is clearly that the optimal choice of  $\tau_{ij}$ depends on the sharp bounds for $|\hat{\omega}_{ij}-\omega_{ij}|$  which are much more difficult to establish than in the covariance matrix case.

Several recent papers considered the estimation of nonparanormal graphical models where the population distribution is non-Gaussian, see Xue and Zou (2012) and Liu, et al. (2012). The nonparanormal model assumes that the variables
follow a joint normal distribution after a set of unknown marginal monotone
transformations. Xue and Zou (2012) estimated the nonparanormal model by applying CLIME (and  graphical lasso, neighborhood Dantzig selector) to the adjusted Spearman's rank correlations. ACLIME can also be used in such a setting. It would be interesting to investigate the properties of the resulting estimator under the  nonparanormal model. Detailed analysis is involved and we leave this as future work.


\section{Proofs}
\label{proof.sec}

In this section we prove the main results, Theorems  \ref{cth1-1} and \ref{cminiupp}, and the key technical results, Lemmas \ref{reduction}, \ref{dffbd} and \ref{affbd}, used in the proof of Theorem \ref{Operlower.bd.thm}. The proof of Lemma \ref{affbd} is involved.
We begin by proving Lemma \ref{var.bnd.lem} stated in Section \ref{method.sec} and collecting  a few additional technical lemmas that will be used in the proofs of the main results.

\subsection{Proof of Lemma \ref{var.bnd.lem} and Additional Technical Lemmas}

\noindent\textbf{Proof of Lemma \ref{var.bnd.lem}} Let $%
\tilde{\Sigma}=(\tilde{\sigma}_{ij})=n^{-1}\sum_{k=1}^{n-1}{X}_{k}%
{X}^{^{\prime}}_{k}$. Note that $\Sigma^{*}$ has the same
distribution as that of $\tilde{\Sigma}$ with ${X}_{k}\sim
N(0,\Sigma) $. So we can replace $\Sigma^{*}$ in Section 2 by $\tilde{%
\Sigma}_{n}=\tilde{\Sigma}+n^{-1}I_{p\times p}$ and assume ${X}%
_{k}\sim N(0,\Sigma)$.
Let
$A_{n}=1-O\left((\log p)^{-1/2}p^{-\delta^{2}/4+1}\right)$ and
set $\lambda _{n}=\delta \sqrt{\log p/n%
}+O((n\log p)^{-1/2})$.
 It suffices to prove that with probability greater than $A_{n}$,
\begin{eqnarray}\label{a6}
\left\vert \sum_{k=1}^{n-1}X_{ki}{X}_{k}^{^{\prime }}\omega _{\cdot
j}\right\vert &\leq &n\lambda _{n}\sqrt{\sigma_{ii}\omega _{jj}%
}\text{, for }i\neq j \cr
\left\vert \sum_{k=1}^{n-1}X_{kj}{X}_{k}^{^{\prime }}\omega _{\cdot
j}-n\right\vert &\leq &n\lambda _{n}\sqrt{\sigma_{jj}\omega
_{jj}}-1\text{, for }i=j.
\end{eqnarray}
Note that $\mathsf{Cov}(\mathbf{X}_{k}^{^{\prime }}\Omega )=\Omega $, $%
\mathsf{Var}({X}_{k}^{^{\prime }}\omega _{\cdot j})=\omega _{jj}$ and
$\mathsf{Cov}(X_{ki}{X}_{k}^{^{\prime }}\omega _{\cdot
j})=\sum_{k=1}^{p}\sigma _{ik}\omega _{kj}=0$ for $i\neq j$. So $X_{ki}$ and
${X}_{k}^{^{\prime }}\omega _{\cdot j}$ are independent. Hence, $\mathbb{E} (X_{ki}{X}_{k}^{^{\prime
}}\omega _{\cdot j})^{3}=0$. By Theorem 5.23 and (5.77) in Petrov (1995),  we
have
\begin{eqnarray}\label{a1}
&&\mathbb{P}\left( \left\vert \sum_{k=1}^{n-1}X_{ki}{X}_{k}^{^{\prime
}}\omega _{\cdot j}\right\vert \geq n\lambda _{n}%
\sqrt{\sigma_{ii}\omega _{jj}}\right) \cr
&&=(1+o(1))\mathbb{P}\left( \left\vert N(0,1)\right\vert \geq \delta \sqrt{\log p}%
\right) \leq C(\log p)^{-1/2}p^{-\delta ^{2}/2}.
\end{eqnarray}%
We next prove the second inequality in (\ref{a6}). We have $\mathbb{E}(X_{kj}%
{X}_{k}^{^{\prime }}\omega _{\cdot j})=1$ and $\mathsf{Var}(X_{kj}%
{X}_{k}^{^{\prime }}\omega _{\cdot j})=\sigma _{jj}\omega _{jj}+1$. Note that $\mathbb{E}\exp (t_{0}(X_{kj}{X}_{k}^{^{\prime }}\omega
_{\cdot j})^{2}/(1+\sigma _{jj}\omega _{jj})\leq c_{0}$ for some absolute
constants $t_{0}$ and $c_{0}$. By Theorem 5.23 in Petrov (1995),
\begin{eqnarray}\label{a2}
\mathbb{P}\Big{(}\Big{|}\sum_{k=1}^{n-1}X_{kj}{X}_{k}^{^{\prime }}\omega _{\cdot
j}-n+1\Big{|}\geq \delta\sqrt{(\sigma_{jj}\omega
_{jj}+1)\log p}\Big{)}\leq C(\log p)^{-1/2}p^{-\delta^{2}/2}.
\end{eqnarray}
Since $1=\mathbb{E}(X_{kj}{X}_{k}^{^{\prime }}\omega _{\cdot j})\leq
\mathbb{E}^{1/2}(X_{kj}{X}_{k}^{^{\prime }}\omega _{\cdot j})^{2}\leq
\sigma _{jj}^{1/2}\omega _{jj}^{1/2}$, we have $\sigma _{jj}\omega _{jj}\geq
1$. This, together with (\ref{a1}) and (\ref{a2}), yields (\ref{a6}).

\begin{lemma}\label{le5}
Let $\hat{\Omega}$ be any estimator of $\Omega$ and set $t_{n}=|\hat{\Omega}-\Omega|_{\infty}$. Then on the event
\begin{eqnarray*}
\left\{|\hat{\omega}_{\cdot j}|_{1}\leq |\omega_{\cdot j}|,\quad%
\mbox{for
$1\leq j\leq p$}\right\},
\end{eqnarray*}
we have
\begin{eqnarray}  \label{a44}
\|\hat{\Omega}-\Omega\|_{1}\leq 12c_{n,p}t_{n}^{1-q}.
\end{eqnarray}
\end{lemma}

\noindent \textbf{Proof.} Define
\begin{equation*}
h_{j}=\hat{\omega}_{\cdot j}-\omega _{\cdot j},\text{ }h_{j}^{1}=(\hat{\omega%
}_{ij}I\{|\hat{\omega}_{ij}|\geq 2t_{n}\};1\leq i\leq p)^{T}-\omega
_{j},~~h_{j}^{2}=h_{j}-h_{j}^{1}.
\end{equation*}%
Then
\begin{equation*}
|\omega _{\cdot j}|_{1}-|h_{j}^{1}|_{1}+|h_{j}^{2}|_{1}\leq |\omega _{\cdot
j}+h_{j}^{1}|_{1}+|h_{j}^{2}|_{1}=|\hat{\omega}_{\cdot j}|_{1}\leq |\omega
_{\cdot j}|_{1},
\end{equation*}%
which implies that $|h_{j}^{2}|_{1}\leq |h_{j}^{1}|_{1}$. This follows that $%
|h_{j}|_{1}\leq 2|h_{j}^{1}|_{1}$. So we only need to upper bound $%
|h_{j}^{1}|_{1}$. We have%
\begin{eqnarray*}
|h_{j}^{1}|_{1} &\leq &\sum_{i=1}^{p}|\hat{\omega}_{ij}-\omega _{ij}|I\{|%
\hat{\omega}_{ij}|\geq 2t_{n}\}+\sum_{i=1}^{p}|\omega _{ij}|I\{|\hat{\omega}%
_{ij}|<2t_{n}\} \\
&\leq &\sum_{i=1}^{p}t_{n}I\{|\omega _{ij}|\geq
t_{n}\}+\sum_{i=1}^{p}|\omega _{ij}|I\{|\omega _{ij}|<3t_{n}\}\leq
4c_{n,p}t_{n}^{1-q}.
\end{eqnarray*}
So (\ref{a44}) holds. \qed

The following Lemma is a classical result. It implies that, if we only
consider estimators of symmetric matrices, an upper bound under the matrix $%
\ell _{1}$ norm is an upper bound for the general matrix $\ell _{w}$ norm
for all $1\leq w\leq \infty $, and a lower bound under the matrix $\ell _{2}$
norm is also a lower bound for the general matrix $\ell _{w}$ norm. We give
a proof to this lemma to be self-contained.

\begin{lemma}
\label{lwl2}Let $A$ be a symmetric matrix, then
\begin{equation*}
\left\Vert A\right\Vert _{2}\leq \left\Vert A\right\Vert _{w}\leq \left\Vert
A\right\Vert _{1}
\end{equation*}%
for all $1\leq w\leq \infty $.
\end{lemma}

\noindent \textbf{Proof of Lemma \ref{lwl2}.} The Riesz-Thorin Interpolation Theorem (See, e.g.,  Thorin, 1948) implies
\begin{equation}
\left\Vert A\right\Vert _{w}\leq \max \left\{ \left\Vert A\right\Vert
_{w_{1}},\left\Vert A\right\Vert _{w_{2}}\right\} \text{,}\ \text{for\ all}\
1\leq w_{1}\leq w\leq w_{2}\leq \infty \text{.}\   \label{RT}
\end{equation}%
Set $w_{1}=1$ and $w_{2}=\infty $, then Equation (\ref{RT}%
) yields $\left\Vert A\right\Vert _{w}\leq \max \left\{ \left\Vert
A\right\Vert _{1},\left\Vert A\right\Vert _{\infty }\right\} $ for all $%
1\leq w\leq \infty $. When $A$ is symmetric, we know $\left\Vert
A\right\Vert _{1}=\left\Vert A\right\Vert _{\infty }$, then immediately we
have $\left\Vert A\right\Vert _{w}\leq \left\Vert A\right\Vert _{1}$. Since %
$2$ is sandwiched between $w$ and $\frac{w}{w-1}$, and $%
\left\Vert A\right\Vert _{w}=\left\Vert A\right\Vert _{\frac{w}{w-1}}$ by
duality, from Equation (\ref{RT}) we have $\left\Vert A\right\Vert _{2}\leq
\left\Vert A\right\Vert _{w}$ for all $1\leq w\leq \infty $ when $A$
symmetric. \qed

\subsection{Proof of Theorems  \ref{cth1-1} and \ref{cminiupp}}

We first prove Theorem \ref{cth1-1}. From Lemma \ref{lwl2} it is enough to consider the $w=1$ case.
 By Lemma 1, we have with probability greater than $A_{n}$,
\begin{eqnarray}\label{a7}
|\hat{\Omega}_{1}-\Omega |_{\infty } &=&|(\Omega \hat{\Sigma}-I_{p\times p})%
\hat{\Omega}_{1}+\Omega (I_{p\times p}-\hat{\Sigma}\hat{\Omega}%
_{1})|_{\infty }\cr
&\leq& C\Vert \hat{\Omega}_{1}\Vert _{1}\sqrt{\frac{\log p%
}{n}}+2\Vert \Omega \Vert _{1}\max_{j}\sigma_{jj}\max_{j}\hat{%
\omega}_{jj}^{1}\sqrt{\frac{\log p}{n}}.
\end{eqnarray}%
We
first assume that $\max_{i}\omega _{ii}>0.5\sqrt{\log p/n}$. By the above
inequality,
\begin{equation*}
\max_{i}\left\vert \frac{\hat{\omega}_{ii}^{1}}{\omega _{ii}}-1\right\vert
\leq Cc_{n,p}\max_{i}\omega _{ii}^{-q}\sqrt{\frac{\log p}{n}}%
+3Mc_{n,p}\max_{i}\omega _{ii}^{-q}\max_{j}\frac{\hat{\omega}_{jj}^{1}}{%
\omega _{jj}}\sqrt{\frac{\log p}{n}}
\end{equation*}%
with probability greater than $A_{n}$. Because $\lambda_{\max}(\Omega)/\lambda_{\min}(\Omega)\leq M_{1}$, we have $\max_{i}\omega _{ii}^{-q}\leq 2(n/\log p)^{q/2}$. Thus by the conditions
in Theorems \ref{cminiupp} and \ref{cth1-1}, we have
\begin{equation*}
\max_{i}\left\vert \frac{\hat{\omega}_{ii}^{1}}{\omega _{ii}}-1\right\vert
=\left\{
\begin{array}{cc}
o(1), & \text{under conditions of Theorem 3} \\
O(1/(\log p)), & \text{under conditions of Theorem 2}%
\end{array}%
\right.
\end{equation*}%
with probability greater than $A_{n}$.
By (\ref{a6}), we can see that, under conditions of Theorem 2, $%
\Omega $ belongs to the feasible set in (\ref{op1}) with probability greater than $A_{n}$. Under conditions of Theorem 3, 
$%
\Omega $ belongs to the feasible set in (\ref{op1}) with probability greater than $1-O\left((\log p)^{-1/2}p^{-\delta^{2}/4+1+o(1)}\right)$. So by a similar
argument as in (\ref{a7}), we can get $|\tilde{\Omega}^{1}-\Omega |_{\infty
}\leq CM_{n,p}\sqrt{\frac{\log p}{n}}$ and $$|\hat{\Omega}-\Omega |_{\infty
}\leq CM_{n,p}\sqrt{\frac{\log p}{n}}.$$ By Lemma \ref{le5}, we see that $$\Vert
\hat{\Omega}-\Omega \Vert _{1}\leq CM_{n,p}^{1-q}c_{n,p}(\log
p/n)^{(1-q)/2}.$$
We consider the case $\max_{i}\omega_{ii}\leq 0.5\sqrt{\log p/n}$. Under this setting, we have $\min_{1\leq i\leq j}\sigma^{*}_{ii}>\sqrt{n/\log p}$ with probability greater than $A_{n}$.  Hence $%
\check{\omega}_{ii}=\sqrt{\log p/n}\geq \omega_{ii}$ and $\Omega$ belongs to
the feasible set in (\ref{op1}) with probability greater than $A_{n}$. So $\|%
\hat{\Omega}\|_{1}\leq \|\Omega\|_{1}\leq Cc_{n,p}(\log
p/n)^{(1-q)/2}$. This proves Theorem \ref{cth1-1}.

To prove Theorem  \ref{cminiupp}, note that $\|\hat{\Omega}\|_{1}\leq \|\tilde{\Omega}^{1}\|_{1}\leq\|\hat{\Sigma}^{-1}\|_{1}\leq np^{1/2}$. We have
\begin{eqnarray*}
\mathbb{E}\left\|\hat{\Omega}-\Omega\right\|^{2}&\leq& CM^{2-2q}_{n,p}c^{2}_{n,p}\left(\frac{\log p}{n}\right)^{1-q}+C(n^{2}p+M^{2-2q}_{n,p}c^{2}_{n,p})p^{-\delta^{2}/4+1+o(1)}(\log p)^{-1/2}\cr
&\leq& CM^{2-2q}_{n,p}c^{2}_{n,p}\left(\frac{\log p}{n}\right)^{1-q}.
\end{eqnarray*}
This proves Theorem \ref{cminiupp}. \qed

\subsection{Proof of Lemma \protect\ref{reduction}}

We first show that the minimax lower bound over all possible estimators is
at the same order of the minimax lower over only estimators of symmetric
matrices under each matrix $\ell _{w}$ norm. For each estimator $\hat{\Omega}
$, we define a projection of $\hat{\Omega}$ to the parameter space $\mathcal{%
F}$,
\begin{equation*}
\hat{\Omega}_{\rm project}=\arg \min_{\Omega \in \mathcal{F}}\left\Vert \hat{%
\Omega}-\Omega \right\Vert _{w},
\end{equation*}%
which is symmetric, then%
\begin{eqnarray}
\sup_{\mathcal{F}}\mathbb{E}\left\Vert \hat{\Omega}_{\rm project}-\Omega
\right\Vert _{w}^{2} &\leq &\sup_{\mathcal{F}}\mathbb{E}\left[ \left\Vert
\hat{\Omega}-\hat{\Omega}_{\rm project}\right\Vert _{w}+\left\Vert \hat{\Omega}%
-\Omega \right\Vert _{w}\right] ^{2}  \notag \\
&\leq &\sup_{\mathcal{F}}\mathbb{E}\left[ \left\Vert \hat{\Omega}-\Omega
\right\Vert _{w}+\left\Vert \hat{\Omega}-\Omega \right\Vert _{w}\right] ^{2}
\notag \\
&=&4\sup_{\mathcal{F}}\mathbb{E}\left\Vert \hat{\Omega}-\Omega \right\Vert
_{w}^{2},  \label{symmest1}
\end{eqnarray}%
where the first inequality follows from the triangle inequality and the
second one follows from the definition of $\hat{\Omega}_{\rm project}$. Since
Equation (\ref{symmest1}) holds for every $\hat{\Omega}$, we have
\begin{equation*}
\inf_{\hat{\Omega},\mbox{\rm symmetric}}\sup_{\mathcal{F}}\mathbb{E}%
\left\Vert \hat{\Omega}-\Omega \right\Vert _{w}^{2}\leq 4\inf_{\hat{\Omega}%
}\sup_{\mathcal{F}}\mathbb{E}\left\Vert \hat{\Omega}-\Omega \right\Vert
_{w}^{2}.
\end{equation*}

From Lemma \ref{lwl2}, we have%
\begin{equation*}
\inf_{\hat{\Omega},\mbox{\rm symmetric}}\sup_{\mathcal{F}}\mathbb{E}%
\left\Vert \hat{\Omega}-\Omega \right\Vert _{w}^{2}\geq \inf_{\hat{\Omega},%
\mbox{\rm symmetric}}\sup_{\mathcal{F}}\mathbb{E}\left\Vert \hat{\Omega}%
-\Omega \right\Vert _{2}^{2}\geq \inf_{\hat{\Omega}}\sup_{\mathcal{F}}%
\mathbb{E}\left\Vert \hat{\Omega}-\Omega \right\Vert _{2}^{2},
\end{equation*}%
which, together with Equation (\ref{symmest1}), establishes Lemma \ref%
{reduction}. \qed

\subsection{Proof of Lemma \protect\ref{dffbd}}

Let $v=\left( v_{i}\right) $ be a column vector with length $p$, and
\begin{equation*}
v_{i}=\left\{
\begin{array}{cc}
1, & p-\lceil p/2\rceil +1\leq i\leq p \\
0, & \text{otherwise}%
\end{array}%
\right.
\end{equation*}%
i.e., $v=\left( 1\left\{ p-\lceil p/2\rceil +1\leq i\leq p\right\} \right)
_{p\times 1}$. Set
\begin{equation*}
w=\left( w_{i}\right) =\left[ \Omega (\theta )-\Omega (\theta ^{\prime })%
\right] v\text{.}
\end{equation*}%
Note that for each $i$, if $\left\vert \gamma _{i}(\theta )-\gamma
_{i}(\theta ^{\prime })\right\vert =1$, we have $\left\vert w_{i}\right\vert
=\frac{M_{n,p}}{2}k\epsilon _{n,p}$. Then there are at least $H(\gamma
(\theta ),\gamma (\theta ^{\prime }))$ number of elements $w_{i}$ with $%
\left\vert w_{i}\right\vert =\frac{M_{n,p}}{2}k\epsilon _{n,p}$, which
implies%
\begin{equation*}
\left\Vert \left[ \Sigma (\theta )-\Sigma (\theta ^{\prime })\right]
v\right\Vert _{2}^{2}\geq H(\gamma (\theta ),\gamma (\theta ^{\prime
}))\cdot \left( \frac{M_{n,p}}{2}k\epsilon _{n,p}\right) ^{2}\text{.}
\end{equation*}%
Since $\left\Vert v\right\Vert ^{2}=\left\lceil p/2\right\rceil \leq p$, the
equation above yields
\begin{equation*}
\left\Vert \Omega (\theta )-\Omega (\theta ^{\prime })\right\Vert ^{2}\geq
\frac{\left\Vert \left[ \Omega (\theta )-\Omega (\theta ^{\prime })\right]
v\right\Vert _{2}^{2}}{\left\Vert v\right\Vert ^{2}}\geq \frac{H(\gamma
(\theta ),\gamma (\theta ^{\prime }))\cdot (\frac{M_{n,p}}{2}k\epsilon
_{n,p})^{2}}{p}\text{,}
\end{equation*}%
i.e.,%
\begin{equation*}
\frac{\left\Vert \Omega (\theta )-\Omega (\theta ^{\prime })\right\Vert ^{2}%
}{H(\gamma (\theta ),\gamma (\theta ^{\prime }))}\geq \frac{%
(M_{n,p}k\epsilon _{n,p})^{2}}{4p}
\end{equation*}%
when $H(\gamma (\theta ),\gamma (\theta ^{\prime }))\geq 1$. \qed


\subsection{Proof of Lemma \protect\ref{affbd}}


Without loss of generality we assume that $M_{n,p}$ is a
constant, since the total variance affinity is scale invariant. The proof of
the bound for the affinity given in Lemma \ref{affbd} is involved. We break
the proof into a few major technical lemmas
Without loss of generality we consider
only the case $i=1$ and prove that there exists a constant $c_{2}>0$ such
that $\left\Vert \mathbb{\bar{P}}_{1,0}\wedge \mathbb{\bar{P}}%
_{1,1}\right\Vert \geq c_{2}$. The following lemma turns the problem of
bounding the total variation affinity into a chi-square distance calculation
on Gaussian mixtures. Define
\begin{equation*}
\Theta _{-1}=\left\{ \left( b,c\right) :\text{there exists a }\theta \in
\Theta \text{ such that }\gamma _{-1}(\theta )=b\text{ and }\lambda
_{-1}(\theta )=c\right\} \text{.}
\end{equation*}%
which is the set of all values of the upper triangular matrix $\Omega \left(
\theta \right) $ could possibly take, with the first row leaving out.

\begin{lemma}
\label{chisquares.lem}If there is \ a constant $c_{2}<1$ such that%
\begin{equation}
\underset{\left( \gamma _{-1},\lambda _{-1}\right) \in \Theta _{-1}}{\mathrm{%
Average}}\left\{ \int \left( \frac{d\mathbb{\bar{P}}_{\left( 1,1,\gamma
_{-1},\lambda _{-1}\right) }}{d\mathbb{\bar{P}}_{\left( 1,0,\gamma
_{-1},\lambda _{-1}\right) }}\right) ^{2}d\mathbb{\bar{P}}_{\left(
1,0,\gamma _{-1},\lambda _{-1}\right) }-1\right\} \leq c_{2}^{2}\text{,}
\label{chi-squarebd}
\end{equation}%
then $\left\Vert \mathbb{\bar{P}}_{1,0}\wedge \mathbb{\bar{P}}%
_{1,1}\right\Vert \geq 1-c_{2}>0$.
\end{lemma}

From the definition of $\mathbb{\bar{P}}_{\left( 1,0,\gamma _{-1},\lambda
_{-1}\right) }$ in Equation (\ref{avepibd}) and $\theta $ in Equation (\ref%
{theta}), $\gamma _{1}=0$ implies $\mathbb{\bar{P}}_{\left( 1,0,\gamma
_{-1},\lambda _{-1}\right) }$ is a single multivariate normal distribution
with a precision matrix,
\begin{equation}
\Omega _{0}=\left(
\begin{array}{cc}
1 & \mathbf{0}_{1\times \left( p-1\right) } \\
\mathbf{0}_{\left( p-1\right) \times 1} & \mathbf{S}_{\left( p-1\right)
\times \left( p-1\right) }%
\end{array}%
\right)  \label{sigma0}
\end{equation}%
where $\mathbf{S}_{\left( p-1\right) \times \left( p-1\right) }=\left(
s_{ij}\right) _{2\leq i,j\leq p}$ is uniquely determined by $\left( \gamma
_{-1},\lambda _{-1}\right) =\left( (\gamma _{2},...,\gamma _{r}),(\lambda
_{2},...,\lambda _{r})\right) $ with%
\begin{equation*}
s_{ij}=\left\{
\begin{array}{cc}
1, & i=j \\
\epsilon _{n,p}, & \gamma _{i}=\text{ }\lambda _{i}\left( j\right) =1 \\
0, & \text{otherwise}%
\end{array}%
\right. \text{.}
\end{equation*}%
Let
\begin{equation*}
\Lambda _{1}\left( c\right) =\left\{ a:\text{there exists a }\theta \in
\Theta \text{ such that }\lambda _{1}(\theta )=a\text{ and }\lambda
_{-1}(\theta )=c\right\} \text{,}
\end{equation*}%
which gives the set of all possible values of the first row with rest of
rows given, i.e., $\lambda _{-1}(\theta )=c$, and define $p_{\lambda _{-1}}=%
\mathrm{Card}\left( \Lambda _{1}\left( \lambda _{-1}\right) \right) $, the
cardinality of all possible $\lambda _{1}$ such that $\left( \lambda
_{1},\lambda _{-1}\right) \in \Lambda $ for the given $\lambda _{-1}$. Then
from definitions in Equations (\ref{avepibd}) and (\ref{theta}) $\mathbb{%
\bar{P}}_{\left( 1,1,\gamma _{-1},\lambda _{-1}\right) }$ is an average of $%
\binom{p_{\lambda _{-1}}}{k}$ multivariate normal distributions with
precision matrices of the following form%
\begin{equation}
\left(
\begin{array}{cc}
1 & \mathbf{r}_{1\times \left( p-1\right) } \\
\mathbf{r}_{\left( p-1\right) \times 1} & \mathbf{S}_{\left( p-1\right)
\times \left( p-1\right) }%
\end{array}%
\right)  \label{matrix.form}
\end{equation}%
where $\left\Vert \mathbf{r}\right\Vert _{0}=k$ with nonzero elements of $r$
equal $\epsilon _{n,p}$ and the submatrix $\mathbf{S}_{(p-1)\times (p-1)}$
is the same as the one for $\Sigma _{0}$ given in (\ref{sigma0}). It is
helpful to observe that $p_{\lambda _{-1}}\geq p/4$. Let $n_{\lambda _{-1}}$
be the number of columns of $\lambda _{-1}$ with column sum equal to $2k$
for which the first row has no choice but to take value $0$ in this column.
Then we have $p_{\lambda _{-1}}=\left\lceil p/2\right\rceil -n_{\lambda
_{-1}}$. Since $n_{\lambda _{-1}}\cdot 2k\leq \left\lceil p/2\right\rceil
\cdot k$, the total number of $1$'s in the upper triangular matrix by the
construction of the parameter set, we thus have $n_{\lambda _{-1}}\leq
\left\lceil p/2\right\rceil /2$, which immediately implies $p_{\lambda
_{-1}}=\left\lceil p/2\right\rceil -n_{\lambda _{-1}}\geq \left\lceil
p/2\right\rceil /2\geq p/4$.

With Lemma \ref{chisquares.lem} in place, it remains to establish Equation (%
\ref{chi-squarebd}) in order to prove Lemma \ref{affbd}. The following lemma
is useful for calculating the cross product terms in the chi-square distance
between Gaussian mixtures. The proof of the lemma is straightforward and is
thus omitted.

\begin{lemma}
\label{crossproduct} Let $g_{i}$ be the density function of $N\left(
0,\Omega _{i}^{-1}\right) $ for $i=0,1$ and $2$. Then%
\begin{equation*}
\int \frac{g_{1}g_{2}}{g_{0}}=\frac{\det \left( I\right) }{\left[ \det
\left( I-\Omega _{1}^{-1}\Omega _{2}^{-1}\left( \Omega _{2}-\Omega
_{0}\right) \left( \Omega _{1}-\Omega _{0}\right) \right) \right] ^{1/2}}%
\text{.}
\end{equation*}
\end{lemma}

Let $\Omega _{i}$, $i=1$ or $2$, be two precision matrices of the form (\ref%
{matrix.form}). Note that $\Omega _{i}$, $i=0,1$ or $2$, differs from each
other only in the first row/column. Then $\Omega _{i}-\Omega _{0}$, $i=1$ or
$2$, has a very simple structure. The nonzero elements only appear in the
first row/column, and in total there are $2k$ nonzero elements. This
property immediately implies the following lemma which makes the problem of
studying the determinant in Lemma \ref{crossproduct} relatively easy.

\begin{lemma}
\label{sigmai} Let $\Omega _{i}$, $i=1$ and $2$, be the precision
matrices of the form (\ref{matrix.form}). Define $J$ to be the number of
overlapping $\epsilon _{n,p}$'s between $\Omega _{1}$ and $\Omega _{2}$ on
the first row, and
\begin{equation*}
Q\overset{\bigtriangleup }{=}\left( q_{ij}\right) _{1\leq i,j\leq p}=\left(
\Omega _{2}-\Omega _{0}\right) \left( \Omega _{1}-\Omega _{0}\right) .
\end{equation*}%
There are index subsets $I_{r}$ and $I_{c}$ in $\left\{ 2,\ldots ,p\right\} $
with $\mathrm{Card}\left( I_{r}\right) =\mathrm{Card}\left( I_{c}\right) =k$
and $\mathrm{Card}\left( I_{r}\cap I_{c}\right) =J$ such that %
\begin{equation*}
q_{ij}=\left\{
\begin{array}{cc}
J\epsilon _{n,p}^{2}, & i=j=1 \\
\epsilon _{n,p}^{2}, & i\in I_{r}\text{ and }j\in I_{c} \\
0, & \text{otherwise}%
\end{array}%
\right.
\end{equation*}%
and the matrix $\left( \Omega _{2}-\Omega _{0}\right) \left( \Omega
_{1}-\Omega _{0}\right) $ has rank $2$ with two identical nonzero
eigenvalues $J\epsilon _{n,p}^{2}$.
\end{lemma}

Let%
\begin{equation}
R_{\lambda _{1},\lambda _{1}^{%
{\acute{}}%
}}^{\gamma _{-1},\lambda _{-1}}=-\log \det \left( I-\Omega _{1}^{-1}\Omega
_{2}^{-1}\left( \Omega _{2}-\Omega _{0}\right) \left( \Omega _{1}-\Omega
_{0}\right) \right),  \label{R}
\end{equation}%
where $\Omega _{0}$ is defined in (\ref{sigma0}) and determined by $\left(
\gamma _{-1},\lambda _{-1}\right) $, and $\Omega _{1}$ and $\Omega _{2}$
have the first row $\lambda _{1}$ and $\lambda _{1}^{%
{\acute{}}%
}$ respectively. We drop the indices $\lambda _{1}$, $\lambda _{1}^{%
{\acute{}}%
}$ and $\left( \gamma _{-1},\lambda _{-1}\right) $ from $\Omega _{i}$ to
simplify the notations. Define
\begin{eqnarray*}
&&\Theta _{-1}\left( a_{1},a_{2}\right) \\
&=&\left\{ 0,1\right\} ^{r-1}\otimes \left\{ c\in \Lambda _{-1}:\text{there
exist }\theta _{i}\in \Theta \text{, }i=1\text{ and }2\text{, such that }%
\lambda _{1}(\theta _{i})=a_{i},\text{ }\lambda _{-1}(\theta _{i})=c\right\}.
\end{eqnarray*}%
It is a subset of $\Theta _{-1}$ in which the element can pick both $a_{1}$
and $a_{2}$ as the first row to form parameters in $\Theta $. From Lemma \ref%
{crossproduct} the left hand side of Equation (\ref{chi-squarebd}) can be
written as
\begin{eqnarray*}
&&\underset{\left( \gamma _{-1},\lambda _{-1}\right) \in \Theta _{-1}}{%
\mathrm{Average}}\left\{ \underset{\lambda _{1},\lambda _{1}^{%
{\acute{}}%
}\in \Lambda _{1}\left( \lambda _{-1}\right) }{\mathrm{Average}}\left[ \exp (%
\frac{n}{2}\cdot R_{\lambda _{1},\lambda _{1}^{%
{\acute{}}%
}}^{\gamma _{-1},\lambda _{-1}})-1\right] \right\} \\
&=&\underset{\lambda _{1},\lambda _{1}^{%
{\acute{}}%
}\in B}{\mathrm{Average}}\left\{ \underset{\left( \gamma _{-1},\lambda
_{-1}\right) \in \Theta _{-1}\left( \lambda _{1},\lambda _{1}^{%
{\acute{}}%
}\right) }{\mathrm{Average}}\left[ \exp (\frac{n}{2}\cdot R_{\lambda
_{1},\lambda _{1}^{%
{\acute{}}%
}}^{\gamma _{-1},\lambda _{-1}})-1\right] \right\} \text{.}
\end{eqnarray*}

The following result shows that $R_{\lambda _{1},\lambda _{1}^{%
{\acute{}}%
}}^{\gamma _{-1},\lambda _{-1}}$ is approximately $-\log \det \left(
I-\left( \Omega _{2}-\Omega _{0}\right) \left( \Omega _{1}-\Omega
_{0}\right) \right) $ which is equal to $-2\log \left( 1-J\epsilon
_{n,p}^{2}\right) $ from Lemma \ref{sigmai}. \ Define%
\begin{equation*}
\Lambda _{1,J}=\left\{ \left( \lambda _{1},\lambda _{1}^{\prime }\right) \in
\Lambda _{1}\otimes \Lambda _{1}:\text{the number of overlapping }\epsilon
_{n,p}\text{'s between }\lambda _{1}\text{and }\lambda _{1}^{\prime }\text{%
is }J\right\} .
\end{equation*}

\begin{lemma}
\label{R.lem}For $R_{\lambda _{1},\lambda _{1}^{%
{\acute{}}%
}}$ defined in Equation (\ref{R}) we have%
\begin{equation}
R_{\lambda _{1},\lambda _{1}^{%
{\acute{}}%
}}^{\gamma _{-1},\lambda _{-1}}=-2\log \left( 1-J\epsilon _{n,p}^{2}\right)
+R_{1,\lambda _{1},\lambda _{1}^{%
{\acute{}}%
}}^{\gamma _{-1},\lambda _{-1}}  \label{Rdecomp}
\end{equation}%
where $R_{1,\lambda _{1},\lambda _{1}^{%
{\acute{}}%
}}^{\gamma _{-1},\lambda _{-1}}$ satisfies
\begin{equation}
\underset{\left( \lambda _{1},\lambda _{1}^{\prime }\right) \in \Lambda
_{1,J}}{\mathrm{Average}}\left[ \underset{\left( \gamma _{-1},\lambda
_{-1}\right) \in \Theta _{-1}\left( \lambda _{1},\lambda _{1}^{%
{\acute{}}%
}\right) }{\mathrm{Average}}\exp \left( \frac{n}{2}R_{1,\lambda _{1},\lambda
_{1}^{%
{\acute{}}%
}}^{\gamma _{-1},\lambda _{-1}}\right) \right] =1+o\left( 1\right) ,
\label{R1.eq}
\end{equation}%
where $J$ is defined in Lemma \ref{sigmai}.
\end{lemma}

\subsubsection{Proof of Equation (\protect\ref{chi-squarebd})}

\label{chisquare.sec}

We are now ready to establish Equation (\ref{chi-squarebd}) which is the key
step in proving Lemma \ref{affbd}. It follows from Equation (\ref{Rdecomp})
in Lemma \ref{R.lem} that
\begin{eqnarray*}
&&\underset{\lambda _{1},\lambda _{1}^{%
{\acute{}}%
}\in B}{\mathrm{Average}}\left\{ \underset{\left( \gamma _{-1},\lambda
_{-1}\right) \in \Theta _{-1}\left( \lambda _{1},\lambda _{1}^{%
{\acute{}}%
}\right) }{\mathrm{Average}}\left[ \exp (\frac{n}{2}R_{\lambda _{1},\lambda
_{1}^{%
{\acute{}}%
}}^{\gamma _{-1},\lambda _{-1}})-1\right] \right\} \\
&=&\underset{J}{\mathrm{Average}}\left\{ \exp \left[ -n\log \left(
1-J\epsilon _{n,p}^{2}\right) \right] \cdot \underset{\left( \lambda
_{1},\lambda _{1}^{\prime }\right) \in \Lambda _{1,J}}{\mathrm{Average}}%
\left[ \underset{\left( \gamma _{-1},\lambda _{-1}\right) \in \Theta
_{-1}\left( \lambda _{1},\lambda _{1}^{%
{\acute{}}%
}\right) }{\mathrm{Average}}\exp (\frac{n}{2}R_{1,\lambda _{1},\lambda _{1}^{%
{\acute{}}%
}}^{\gamma _{-1},\lambda _{-1}})\right] -1\right\} \text{.}
\end{eqnarray*}%
Recall that $J$ is the number of overlapping $\epsilon _{n,p}$'s between $%
\Sigma _{1}$ and $\Sigma _{2}$ on the first row. It can be shown that $J$
has the hypergeometric distribution with
\begin{equation}
\mathbb{P}\left( \text{number of overlapping }\epsilon _{n,p}\text{'s}%
=j\right) =\binom{k}{j}\binom{p_{\lambda _{-1}}-k}{k-j}/\binom{p_{\lambda
_{-1}}}{k}\leq \left( \frac{k^{2}}{p_{\lambda _{-1}}-k}\right) ^{j}\text{.}
\label{plambda}
\end{equation}
Equation (\ref{plambda}) and Lemma \ref{R.lem}, togetehr with Equation (\ref%
{epsilonnp}), imply%
\begin{eqnarray}
&&\underset{\left( \gamma _{-1},\lambda _{-1}\right) \in \Theta _{-1}}{%
\mathrm{Average}}\left\{ \int \left( \frac{d\mathbb{\bar{P}}_{\left(
1,1,\gamma _{-1},\lambda _{-1}\right) }}{d\mathbb{\bar{P}}_{\left(
1,0,\gamma _{-1},\lambda _{-1}\right) }}\right) ^{2}d\mathbb{\bar{P}}%
_{\left( 1,0,\gamma _{-1},\lambda _{-1}\right) }-1\right\}  \notag \\
&\leq &\sum_{j\geq 0}\left( \frac{k^{2}}{p/4-1-k}\right) ^{j}\left\{ \exp %
\left[ -n\log \left( 1-j\epsilon _{n,p}^{2}\right) \right] \cdot \left(
1+o\left( 1\right) \right) -1\right\}  \label{limitk} \\
&\leq &\left( 1+o\left( 1\right) \right) \sum_{j\geq 1}\left( p^{\frac{\beta
-1}{\beta }}\right) ^{-j}\exp \left[ 2j\left( \upsilon ^{2}\log p\right) %
\right] +o\left( 1\right)  \notag \\
&\leq &C\sum_{j\geq 1}\left( p^{\frac{\beta -1}{\beta }-2\upsilon
^{2}}\right) ^{-j}+o\left( 1\right) <c_{2}^{2},  \notag
\end{eqnarray}%
where the last step follows from $\upsilon ^{2}<\frac{\beta -1}{8\beta },$
and $k^{2}\leq \left[ c_{n,p}(M_{n,p}\epsilon _{n,p})^{-q}\right]
^{2}=O\left( \frac{n}{\log ^{3}p}\right) =o\left( \frac{p^{1/\beta }}{\log p}%
\right) $ from Equations (\ref{kdef}) and (\ref{cnp}) and the condition $p>c_{1}n^{\beta }$ for some $\beta >1$%
, and $c_{2}$ is a positive constant. \quad
\hbox{\vrule width
4pt height 6pt depth 1.5pt}

\subsection{Proof of Lemma \protect\ref{R.lem}}

\label{tech.lemmas.sec}

Let%
\begin{equation}
A=\left( I-\Omega _{1}^{-1}\Omega _{2}^{-1}\right) \left( \Omega _{2}-\Omega
_{0}\right) \left( \Omega _{1}-\Omega _{0}\right) \left[ I-\left( \Omega
_{2}-\Omega _{0}\right) \left( \Omega _{1}-\Omega _{0}\right) \right] ^{-1}%
\text{.}  \label{Adef}
\end{equation}%
Since $\left\Vert \Omega _{i}-\Omega _{0}\right\Vert \leq \left\Vert \Omega
_{i}-\Omega _{0}\right\Vert _{1}\leq 2k\epsilon _{n,p}=o\left( 1/\log
p\right) $ from Equation(\ref{kepsilon}), it is easy to see that%
\begin{equation}
\left\Vert A\right\Vert =O\left( k\epsilon _{n,p}\right) =o\left( 1\right)
\text{.}  \label{Aeig}
\end{equation}%
Define%
\begin{equation*}
R_{1,\lambda _{1},\lambda _{1}^{%
{\acute{}}%
}}^{\gamma _{-1},\lambda _{-1}}=-\log \det \left( I-A\right) \text{.}
\end{equation*}%
Then we can rewrite $R_{\lambda _{1},\lambda _{1}^{%
{\acute{}}%
}}^{\gamma _{-1},\lambda _{-1}}$ as follows%
\begin{eqnarray}
R_{\lambda _{1},\lambda _{1}^{%
{\acute{}}%
}}^{\gamma _{-1},\lambda _{-1}} &=&-\log \det \left( I-\Omega
_{1}^{-1}\Omega _{2}^{-1}\left( \Omega _{2}-\Omega _{0}\right) \left( \Omega
_{1}-\Omega _{0}\right) \right)  \notag \\
&=&-\log \det \left( \left[ I-A\right] \cdot \left[ I-\left( \Omega
_{2}-\Omega _{0}\right) \left( \Omega _{1}-\Omega _{0}\right) \right] \right)
\notag \\
&=&-\log \det \left[ I-\left( \Omega _{2}-\Omega _{0}\right) \left( \Omega
_{1}-\Omega _{0}\right) \right] -\log \det \left( I-A\right)  \notag \\
&=&-2\log \left( 1-J\epsilon _{n,p}^{2}\right) +R_{1,\lambda _{1},\lambda
_{1}^{%
{\acute{}}%
}}^{\gamma _{-1},\lambda _{-1}}\text{,}  \label{minsR}
\end{eqnarray}%
where the last equation follows from Lemma \ref{sigmai}. To establish Lemma %
\ref{R.lem} it is enough to establish Equation (\ref{R1.eq}).

Let
\begin{eqnarray*}
B_{1} &=&\left( b_{1,ij}\right) _{p\times p}=\left[ I-\left( \Omega
_{1}-\Omega _{0}+I\right) ^{-1}\left( \Omega _{2}-\Omega _{0}+I\right) ^{-1}%
\right] \text{, } \\
B_{2} &=&\left( b_{2,ij}\right) _{p\times p}=\left( \Omega _{2}-\Omega
_{0}\right) \left( \Omega _{1}-\Omega _{0}\right) \left[ I-\left( \Omega
_{2}-\Omega _{0}\right) \left( \Omega _{1}-\Omega _{0}\right) \right] ^{-1},
\end{eqnarray*}%
and define%
\begin{equation}
A_{1}=\left( a_{1,ij}\right) =B_{1}B_{2}\text{.}  \label{A1eig}
\end{equation}%
Similar to Equation (\ref{Aeig}), we have $\left\Vert A_{1}\right\Vert
=o\left( 1\right) $, and write%
\begin{equation*}
R_{1,\lambda _{1},\lambda _{1}^{%
{\acute{}}%
}}^{\gamma _{-1},\lambda _{-1}}=-\log \det \left( I-A_{1}\right) -\log \det %
\left[ \left( I-A_{1}\right) ^{-1}\left( I-A\right) \right] \text{.}
\end{equation*}%
To establish Equation (\ref{R1.eq}), it is enough to show that
\begin{equation}
\exp \left[ -\frac{n}{2}\log \det \left( I-A_{1}\right) \right] =1+o\left(
1\right) ,  \label{A1bd}
\end{equation}%
and that
\begin{equation}
\underset{\left( \gamma _{-1},\lambda _{-1}\right) \in \Theta _{-1}\left(
\lambda _{1},\lambda _{1}^{%
{\acute{}}%
}\right) }{\mathrm{Average}}\exp \left( -\frac{n}{2}\log \det \left[ \left(
I-A_{1}\right) ^{-1}\left( I-A\right) \right] \right) =1+o\left( 1\right) .
\label{A1Abd}
\end{equation}

The proof for Equation (\ref{A1bd}) is as follows. Write%
\begin{equation*}
\Omega _{1}-\Omega _{0}=\left(
\begin{array}{cc}
0 & \mathbf{v}_{1\times \left( p-1\right) } \\
\left( \mathbf{v}_{1\times \left( p-1\right) }\right) ^{T} & \mathbf{0}%
_{\left( p-1\right) \times \left( p-1\right) }%
\end{array}%
\right) \text{, and }\Omega _{2}-\Omega _{0}=\left(
\begin{array}{cc}
0 & \mathbf{v}_{1\times \left( p-1\right) }^{\ast } \\
\left( \mathbf{v}_{1\times \left( p-1\right) }^{\ast }\right) ^{T} & \mathbf{%
0}_{\left( p-1\right) \times \left( p-1\right) }%
\end{array}%
\right)
\end{equation*}%
where $\mathbf{v}_{1\times \left( p-1\right) }=\left( v_{j}\right) _{2\leq
j\leq p}$ satisfies $v_{j}=0$ for $2\leq j\leq p-r$ and $v_{j}=0$ or $1$ for
$p-r+1\leq j\leq p$ with $\left\Vert \mathbf{v}\right\Vert _{0}=k$, and $%
\mathbf{v}_{1\times \left( p-1\right) }^{\ast }=\left( v_{j}^{\ast }\right)
_{2\leq j\leq p}$ satisfies a similar property. Without loss of generality
we consider only a special case with
\begin{equation*}
v_{j}=\left\{
\begin{tabular}{ll}
$1,$ & $p-r+1\leq j\leq p-r+k$ \\
$0,$ & otherwise%
\end{tabular}%
\right. \text{, and }v_{j}^{\ast }=\left\{
\begin{tabular}{ll}
$1,$ & $p-r+k-J\leq j\leq p-r+2k-J$ \\
$0,$ & otherwise%
\end{tabular}%
\right. .
\end{equation*}%
Note that $B_{1}$ can be written as a polynomial of $\Omega _{1}-\Omega _{0}$
and $\Omega _{2}-\Omega _{0}$, and $B_{2}$ can be written as a polynomial of
$\left( \Omega _{2}-\Omega _{0}\right) \left( \Omega _{1}-\Omega _{0}\right)
.$ By a straightforward calculation it can be shown that%
\begin{equation*}
\left\vert b_{1,ij}\right\vert =\left\{
\begin{array}{cc}
O\left( \epsilon _{n,p}\right) , &
\begin{array}{c}
i=1\text{ and }p-r+1\leq j\leq p-r+2k-J\text{, } \\
\text{or }j=1\text{ and }p-r+1\leq i\leq p-r+2k-J\text{, or }i=j=1%
\end{array}
\\
O\left( \epsilon _{n,p}^{2}\right) , & \text{ }p-r+1\leq i\leq p-r+2k-J\text{%
, and }p-r+1\leq j\leq p-r+2k-J \\
0, & \text{otherwise}%
\end{array}%
\right. ,
\end{equation*}%
and
\begin{equation*}
0\leq b_{2,ij}=\left\{
\begin{array}{cc}
O\left( J\epsilon _{n,p}^{2}\right) , & i=j=1 \\
\tau _{n,p}, & \text{ }p-r+1\leq i\leq p-r+k\text{, and }p-r+k-J\leq j-1\leq
p-r+2k-J \\
0, & \text{otherwise}%
\end{array}%
\right. \text{,}
\end{equation*}%
where $\tau _{n,p}=O\left( \epsilon _{n,p}^{2}\right) $, which implies%
\begin{equation*}
\left\vert a_{1,ij}\right\vert =\left\{
\begin{tabular}{ll}
$O\left( k\epsilon _{n,p}^{3}\right) ,$ & $i=1\text{ and }p-r+k-J\leq
j-1\leq p-r+2k-J\text{, or }i=j=1\text{ }$ \\
$O\left( J\epsilon _{n,p}^{3}\right) ,$ & $j=1\text{ and }p-r+1\leq i\leq
p-r+2k-J$ \\
$O\left( k\epsilon _{n,p}^{4}\right) ,$ & $\text{ }p-r+1\leq i\leq p-r+2k-J%
\text{, and }p-r+k-J\leq j-1\leq p-r+2k-J$ \\
$0,$ & $\text{otherwise}$%
\end{tabular}%
\right. \text{.}
\end{equation*}%
Note that $\mathrm{rank}\left( A_{1}\right) \leq 2$ due to the simple
structure of $\left( \Omega _{2}-\Omega _{0}\right) \left( \Omega
_{1}-\Omega _{0}\right) $. Let $A_{2}=\left( a_{2,ij}\right) $ with%
\begin{equation*}
\left\vert a_{2,ij}\right\vert =\left\{
\begin{tabular}{ll}
$O\left( k\epsilon _{n,p}^{3}\right) ,$ & $i=1\text{ and }$ $j=1$ \\
$O\left( k\epsilon _{n,p}^{4}+Jk\epsilon _{n,p}^{6}\right) ,$ & $\text{ }%
p-r+1\leq i\leq p-r+2k-J$\\
&$\text{and }p-r+k-J\leq j-1\leq p-r+2k-J$ \\
$0,$ & $\text{otherwise}$%
\end{tabular}%
\right. \text{,}
\end{equation*}%
and $\mathrm{rank}\left( A_{2}\right) \leq 4$ by eliminating the non-zero
off-diagonal elements of the first row and column of $A_{1}$, and
\begin{equation*}
\exp \left[ -\frac{n}{2}\log \det \left( I-A_{1}\right) \right] =\exp \left[
-\frac{n}{2}\log \det \left( I-A_{2}\right) \right] \text{.}
\end{equation*}
We can show that all eigenvalues of $A_{1}^{\ast }$ are $O\left(
Jk^{2}\epsilon _{n,p}^{6}+k^{2}\epsilon _{n,p}^{4}+k\epsilon
_{n,p}^{3}\right) $. Since $k\epsilon _{n,p}=o\left( 1/\log p\right) ,$ then
\begin{equation*}
nk\epsilon _{n,p}^{3}=\upsilon ^{3}\frac{k\left( \log p\right) ^{1/2}}{\sqrt{%
n}}\log p=o\left( 1\right)
\end{equation*}%
which implies%
\begin{equation*}
n\left( Jk^{2}\epsilon _{n,p}^{6}+k^{2}\epsilon _{n,p}^{4}+k\epsilon
_{n,p}^{3}\right) =o\left( 1\right) \text{.}
\end{equation*}%
Thus%
\begin{equation*}
\exp \left[ -\frac{n}{2}\log \det \left( I-A_{1}\right) \right] =1+o\left(
1\right) \text{.}
\end{equation*}

Now we establish Equation (\ref{A1Abd}), which, together with Equation (\ref%
{A1bd}), yields Equation (\ref{R1.eq}) and thus Lemma \ref{R.lem} is
established. Write%
\begin{eqnarray*}
&&\left( I-A_{1}\right) ^{-1}\left( I-A\right) -I=\left( I-A_{1}\right)
^{-1} \left[ \left( I-A\right) -\left( I-A_{1}\right) \right] =\left(
I-A_{1}\right) ^{-1}\left( A_{1}-A\right) \\
&=&\left( I-A_{1}\right) ^{-1}\left[ \Omega _{1}^{-1}\Omega _{2}^{-1}-\left(
\Omega _{1}-\Omega _{0}+I\right) ^{-1}\left( \Omega _{2}-\Omega
_{0}+I\right) ^{-1}\right] \\
&&\cdot \left( \Omega _{2}-\Omega _{0}\right) \left( \Omega _{1}-\Omega
_{0}\right) \left[ I-\left( \Omega _{2}-\Omega _{0}\right) \left( \Omega
_{1}-\Omega _{0}\right) \right] ^{-1}
\end{eqnarray*}%
where%
\begin{eqnarray*}
&&\Omega _{1}^{-1}\Omega _{2}^{-1}-\left( \Omega _{1}-\Omega _{0}+I\right)
^{-1}\left( \Omega _{2}-\Omega _{0}+I\right) ^{-1} \\
&=&\Omega _{1}^{-1}\Omega _{2}^{-1}\left[ \left( \Omega _{2}-\Omega
_{0}+I\right) \left( \Omega _{1}-\Omega _{0}+I\right) -\Omega _{2}\Omega _{1}%
\right] \left( \Omega _{1}-\Omega _{0}+I\right) ^{-1}\left( \Omega
_{2}-\Omega _{0}+I\right) ^{-1} \\
&=&\Omega _{1}^{-1}\Omega _{2}^{-1}\left[ \left( -\Omega _{0}+I\right)
\Omega _{1}+\Omega _{2}\left( -\Omega _{0}+I\right) +\left( -\Omega
_{0}+I\right) ^{2}\right] \left( \Omega _{1}-\Omega _{0}+I\right)
^{-1}\left( \Omega _{2}-\Omega _{0}+I\right) ^{-1}\text{.}
\end{eqnarray*}%
It is important to observe that $\mathrm{rank}\left( \left( I-A_{1}\right)
^{-1}\left( I-A\right) -I\right) \leq 2$ again due to the simple structure
of $\left( \Omega _{2}-\Omega _{0}\right) \left( \Omega _{1}-\Omega
_{0}\right) $, then $-\log \det \left[ \left( I-A_{1}\right) ^{-1}\left(
I-A\right) \right] $ is determined by at most two nonzero eigenvalues, which
are bounded by%
\begin{equation}
\left\Vert \left( I-A_{1}\right) ^{-1}\left( I-A\right) -I\right\Vert
=\left( 1+o\left( 1\right) \right) \left\Vert \left( I-\Omega _{0}\right)
\left( \Omega _{2}-\Omega _{0}\right) \left( \Omega _{1}-\Omega _{0}\right)
\right\Vert \text{.}  \label{AStartbd}
\end{equation}%
Note that $\left\Vert \left( I-A_{1}\right) ^{-1}\left( I-A\right)
-I\right\Vert =o\left( 1\right) $, and%
\begin{equation*}
\left\vert \log \left( 1-x\right) \right\vert \leq 2\left\vert x\right\vert
\text{, for }\left\vert x\right\vert <1/3,
\end{equation*}%
which implies%
\begin{equation*}
\left\vert -\log \det \left[ \left( I-A_{1}\right) ^{-1}\left( I-A\right) %
\right] \right\vert \leq 2\left\Vert \left( I-A_{1}\right) ^{-1}\left(
I-A\right) -I\right\Vert \text{,}
\end{equation*}%
i.e.,%
\begin{equation*}
\exp \left( \frac{n}{2}\cdot -\log \det \left[ \left( I-A_{1}\right)
^{-1}\left( I-A\right) \right] \right) \leq \exp \left( n\left\Vert \left(
I-A_{1}\right) ^{-1}\left( I-A\right) -I\right\Vert \right) \text{.}
\end{equation*}%
Define%
\begin{equation*}
A_{\ast }=\left( I-\Omega _{0}\right) \left( \Omega _{2}-\Omega _{0}\right)
\left( \Omega _{1}-\Omega _{0}\right) ,
\end{equation*}%
then%
\begin{equation*}
\exp \left( \frac{n}{2}\cdot -\log \det \left[ \left( I-A_{1}\right)
^{-1}\left( I-A\right) \right] \right) \leq \exp \left( (1+o(1))n\left\Vert A_{\ast }\right\Vert \right)
\end{equation*}%
from Equations (\ref{AStartbd}). It is then sufficient to show

\begin{equation*}
\underset{\left( \gamma _{-1},\lambda _{-1}\right) \in \Theta _{-1}\left(
\lambda _{1},\lambda _{1}^{%
{\acute{}}%
}\right) }{\mathrm{Average}}\exp \left( 2n\left\Vert A_{\ast }\right\Vert
\right) =1+o\left( 1\right)
\end{equation*}%
where $\left\Vert A_{\ast }\right\Vert $ depends on the values of $\lambda
_{1},\lambda _{1}^{%
{\acute{}}%
}$ and $\left( \gamma _{-1},\lambda _{-1}\right) $. We dropped the indices $%
\lambda _{1}$, $\lambda _{1}^{%
{\acute{}}%
}$ and $\left( \gamma _{-1},\lambda _{-1}\right) $ from $A$ to simplify the
notations.

Let $E_{m}=\left\{ 1,2,\ldots ,r\right\} /\left\{ 1,m\right\} $. Let $%
n_{\lambda _{E_{m}}}$ be the number of columns of $\lambda _{E_{m}}$ with
column sum at least $2k-2$ for which two rows can not freely take value $0$
or $1$ in this column. Then we have $p_{\lambda _{E_{m}}}=\left\lceil
p/2\right\rceil -n_{\lambda _{E_{m}}}$. Without loss of generality we assume
that $k\geq 3$. Since $n_{\lambda _{E_{m}}}\cdot \left( 2k-2\right) \leq
\left\lceil p/2\right\rceil \cdot k$, the total number of $1$'s in the upper
triangular matrix by the construction of the parameter set, we thus have $%
n_{\lambda _{E_{m}}}\leq \left\lceil p/2\right\rceil \cdot \frac{3}{4}$,
which immediately implies $p_{\lambda _{E_{m}}}=\left\lceil p/2\right\rceil
-n_{\lambda _{E_{m}}}\geq \left\lceil p/2\right\rceil \frac{1}{4}\geq p/8$.
Thus we have%
\begin{eqnarray*}
\mathbb{P}\left( \left\Vert A_{\ast }\right\Vert \geq 2t\cdot \epsilon
_{n,p}\cdot k\epsilon _{n,p}^{2}\right) &\leq &\mathbb{P}\left( \left\Vert
A_{\ast }\right\Vert _{1}\geq 2t\cdot \epsilon _{n,p}\cdot k\epsilon
_{n,p}^{2}\right) \\
&\leq &\sum_{m}\underset{\lambda _{E_{m}}}{\mathrm{Average}}\frac{\binom{k}{t%
}\binom{p_{\lambda _{E_{m}}}}{k-t}}{\binom{p_{\lambda _{E_{m}}}}{k}}\leq
p\left( \frac{k^{2}}{p/8-k}\right) ^{t}
\end{eqnarray*}%
from Equation (\ref{plambda}), which immediately implies%
\begin{eqnarray*}
&&\underset{\left( \gamma _{-1},\lambda _{-1}\right) \in \Theta _{-1}\left(
\lambda _{1},\lambda _{1}^{%
{\acute{}}%
}\right) }{\mathrm{Average}}\exp \left( 2n\left\Vert A_{\ast }\right\Vert
\right) \\
&\leq &\exp \left( 4n\cdot \frac{2\left( \beta -1\right) }{\beta }\cdot
\epsilon _{n,p}\cdot k\epsilon _{n,p}^{2}\right) +\int_{\frac{2\left( \beta
-1\right) }{\beta }}^{\infty }\exp \left( 2n\cdot 2t\cdot \epsilon
_{n,p}\cdot k\epsilon _{n,p}^{2}\right) p\left( \frac{k^{2}}{p/8-k}\right)
^{t}dt \\
&=&\exp \left( \frac{8\left( \beta -1\right) }{\beta }nk\epsilon
_{n,p}^{3}\right) +\int_{\frac{2\left( \beta -1\right) }{\beta }}^{\infty
}\exp \left[ \log p+t\left( 4nk\epsilon _{n,p}^{3}-\log \frac{k^{2}}{p/8-k}%
\right) \right] dt \\
&=&1+o\left( 1\right) \text{,}
\end{eqnarray*}%
where the last step is an immediate consequence of the following two
equations,%
\begin{equation*}
nk\epsilon _{n,p}^{3}=o\left( 1\right)
\end{equation*}%
and
\begin{equation*}
\left( 1+o\left( 1\right) \right) 2\log p\leq t\log \frac{p/8-1-k}{k^{2}}%
\text{ , for }t\geq \frac{2\left( \beta -1\right) }{\beta }
\end{equation*}%
which follow from $k^{2}=O\left( n\right) =O\left( p^{1/\beta }\right) $ from Equation (\ref{kdef}) and the condition $p>c_{1}n^{\beta }$ for some $\beta >1$. \qed


\end{document}